\documentclass[1pt]{elsarticle}
\usepackage{amsmath,amsfonts, amssymb}

\newcommand{\C} {\mathbb{C}^N}
\newcommand{\CP} {\mathbb{C}P^{N-1}}
\newcommand{\be} {\begin {equation}}
\newcommand {\ee} {\end {equation}}
\newcommand{\ba} {\begin {array}}
\newcommand {\ea} {\end {array}}
\newcommand{\bea}{\begin{eqnarray}}
\newcommand{\eea}{\end{eqnarray}}
\newcommand{\p}{\partial_+}
\newcommand{\bp}{\partial_-}
\newcommand{\nn}{\nonumber \\}
\newcommand{\0}{{\bf \varnothing} }

\newcommand{\Pb}{{\bf P}}
\newcommand{\X}{{\bf X} }
\newcommand{\Y}{{\bf Y} }
\newcommand{\N}{{\bf N} }
\renewcommand{\H}{{\bf H} }
\newcommand{\F}{{\bf F} }
\renewcommand{\v}{{v}}
\newcommand{\w}{{w}}
\newcommand{\x}{{x}}
\newcommand{\y}{{y}}
\newcommand{\z}{{z}}
\newcommand{\e}{{e}}
\newcommand{\f}{{f}}
\newcommand{\QED}{\qed\\}
\newtheorem{theorem}{Theorem}
\newtheorem{definition}{Definition}
\newtheorem{property}{Property}
\newtheorem{lemma}{Lemma}

\renewcommand{\vec}{{}}

\journal{ArXiv}
\begin{document}
\begin{frontmatter}
\title{Analysis of $\mathbb{C}P^{N-1}$ sigma models via projective structure}

 \author[label1]{S. Post}
 \author[label1,label2]{A. M. Grundland}
  \address[label1]{Centre de Recherches Math\'ematiques, Universit\'e de Montreal. \\C.P. 6128 succ. Centre-Ville, Montreal (QC) H3C 3J7, Canada.}
\address[label2]{ Department of Mathematics and Computer Sciences,\\ Universit\'e de Qu\'ebec, Trois-Riviers. CP500 (QC) G9A 5H7, Canada. }

\begin{abstract}
In this paper, we study rank-1 projector solutions to the completely integrable Euclidean $\CP$ sigma model in two dimension and their associated surfaces immersed in the $su(N)$ Lie algebra. 
We reinterpret and generalize the proof of A.M. Din and W.J. Zakzrewski \cite{DinZak1980Gen} that any  solution for the $\CP$ sigma model defined on the Riemann sphere with finite action can be written as a raising operator acting on a holomorphic one, or a lowering operator acting on a antiholomorphic one. Our proof is formulated in terms  of rank-1 Hermitian projectors so it is explicitly gauge invariant and gives new results on the structure of the corresponding sequence of rank-1 projectors.  Next, we analyze surfaces associated with the $\CP$ models defined using the Generalized Weierstrass Formula for immersion, introduced by B. Konopelchenko \cite{Kono1996}. We show that the surfaces are conformally parameterized by the Lagrangian density with finite area equal to the action of the model and express several other geometrical characteristics of the surface in terms of the Lagrangian density and topological charge density of the model. We demonstrate that any such surface must be orthogonal to  the sequence of projectors defined by repeated application of the raising and lowering operators.  Finally, we provide necessary and sufficient conditions that a surface be related to a $\CP$ sigma model. 
\end{abstract}


\end{frontmatter}

\section{Introduction}
The solutions of the $\CP$ sigma model defined on the extended complex plane with finite action can be characterized by a holomorphic vector and a finite sequence of vectors built using raising operators. This fact was first discover by  A.M. Din and W. Zakzrewski \cite{DinZak1980Gen, DinZak1980prop},  H.J. Borchers and W.D. Garber \cite{BorchersGarber},  next by R. Sasaki \cite{ Sasaki1983}, and later discussed by J. Eells and J. Wood  \cite{EellsWood1983}. It was shown by A.V. Mikhailov and V.E. Zakharov that  the Euler-Lagrange (E-L) equation for vectors in homogenous coordinates can be written in terms of projectors to obtain a new form of the Euler-Lagrange  equations  and they introduced a linear spectral problem whose compatibility  corresponds precisely with such equations \cite{Mik1986, ZakMik1979}. From a complete background on the $\CP$ sigma model and history (see e.g. \cite{Bobbook, ChernWolf1987, Guestbook, Heleinbook2001,  HeleinWoodbook,  ManSutbook, Polybook, Uhlen1989, Wardbook,  Zakbook}).

In this paper, we complete the theory by showing that any projector solution to the (projector) Euler-Lagrange  equations can be given in such a way. More precisely, we show that all rank-1 Hermetian projector solution to the  $\CP$ sigma model defined on the  with finite action are obtained through successive application of a creation (annihilation) operator acting on projectors which map into a direction of a  holomorphic (anti-holomorphic) vector. Though some of the results are already known and in use in current literature, citations above and  \cite{GoldGrund2010}, this reconstructed proof completes the theory and has the advantage that it is written completely in terms of projectors, allowing for full utilization of the projective structures of $\CP$ in a consistent way. This proof not only gives necessary and sufficient conditions for projectors to be in a certain class of solutions for the $\CP$ sigma model but also allows further analysis of the model using the new results on the differential structures of projectors developed for use in the proof. Consequently, these results allow us to explore geometrical characteristics of 2D surfaces immersed in the $su(N)$ algebra and express them in terms of the projectors. Such geometric properties include the Gauss curvature, the mean curvature vectors and some global characteristics of 2D surfaces such as the Willmore functional and the Euler-Poincare characters. Moreover, the results also allow us to express quantities associated with the sigma model, such as the Lagrangian density, action and topological charge, in terms of geometrical properties of the surfaces, such as surface area, Willmore functional, Euler-Poincare etc. 

In studying sigma models on Euclidean space, one is interested in maps $\z:\Omega \subset \mathbb{R}^2\rightarrow \C$ with $\z^\dagger \z=1$ which are stationary points of the action functional \cite{Zakbook}
\be \label{Sz} S(\z)=\sum_{\mu =1,2} \int_{\mathbb{R}^2}(D_\mu \z)^\dagger D_\mu \z \ d\xi_{1}d\xi_2,\ee
given by solutions to the corresponding Euler-Lagrange equations 
$$\sum_{\mu =1,2} D_\mu D_\mu \z+\z((D_\mu \z)^\dagger D_\mu \z)=0 $$
where $D_\mu z$ defined according to the formula
$$D_\mu \z=\partial_\mu \z-(\z^\dagger \partial_\mu \z)\z, \qquad \partial_\mu=\frac{\partial}{\partial \xi_\mu}.$$ 
Notice that, for any $k :\Omega \rightarrow \mathbb{C},$ $D_\mu(k\z)= k D_\mu(\z)$ and so  the action functional defined by \eqref{Sz} is invariant under the local $U(1)$ transformation induced by $k$ with $|k|=1.$ Furthermore, we can expand the model to act on any non-zero vector $\f:\Omega\rightarrow \C \diagdown \lbrace \0 \rbrace $ by 
$$D_\mu \f=\partial_\mu \f-\frac{\f^\dagger \partial_\mu \f}{\f^\dagger \f}\f,$$ 
so that the action becomes
\be \label{Sf} S(\f)= \sum_{\mu =1,2}  \int_{\mathbb{R}^2}\frac{(D_\mu \f)^\dagger D_\mu \f}{\f^\dagger \f}d\xi_1d\xi_2,\qquad  (\xi_1,\xi_2)\in \mathbb{R}^2.\ee
We see immediately that  equations \eqref{Sz} and \eqref{Sf} are consistent with $\z=\f/|\f |, \ |\f|\equiv(\f^\dagger \f)^{1/2}$. Also, for any $\f,\vec{g}:\Omega \rightarrow \C\diagdown \lbrace \0 \rbrace $ where $\f=k\vec{g}$ for some $k:\Omega\rightarrow \mathbb{C}\diagdown \lbrace 0 \rbrace$ the actions are the same,  $S(\f)=S(k\f)=S(\vec{g}),$ and  so our model is defined on equivelence class of $\CP.$ Finally, we note that the action $S(\f)$ is often given in the literature, \cite{Zakbook}, as having values $2\pi \alpha_0$ where $\alpha_0$ is the degree of the vector field $\f.$

In this paper, we provide an analysis of these models by exploiting the projective structures of $\CP$ though the projector formalism, ie in terms of $\Pb =(\f\otimes \f^\dagger)/(\f^\dagger \f).$ The language of projectors ensure invariance under multiplication of $f$ by any non-zero scalar function and freedom of singularites which could occur in the description in terms of unnormalized vector fields $f$ \cite{GoldGrund2010}. Moreover, ther previously known properties of the $\CP$ sigma model and the associated surfaces take a compact and simple form when expressed in terms of projectors. 

The paper is organized as follows. 
In Section 2, we introduce a new sesquelinear product, semi-norm and covariant derivative to act on the space of projectors and proven new results  associated with them such as the Cauchy-Schwarz inequality and  the triangle inequality for the semi-norm. These structures have been defined in such a way that the covariant derivative is compatible with the sesquilinear product while simultaneously being orthogonal to the original operator. We also apply the concept of (anti-)holomorphicity to projectors and present a new analysis of structures, particularly in relation to a set of raising and lowering operators $\Pi_\pm, $ as introduced in \cite{GoldGrund2010}. From example, we prove that a rank-1 Hermitian projector will project onto a direction in $\CP$ which has a holomorphic representative if and only if the lowering operator $\Pi_-$ annihilates the projector. 

In Section 3, we provide the basic background on $\CP$ sigma models and write the physical quanities and Euler-Lagrange equations explicitly in terms of the covariant derivatives defined in Section 2. 

In Section 4, we start with a a rank-1 Hermetian projector which is a solution to the Euler-Lagrange  equations defined on the extended complex plane with finite action and create via raising and lowering operators a finite set of rank-1 Hermitian projectors which are solutions to the Euler-Lagrange  equations with finite action and are mutually orthogonal.  More specifically, in Theorem 4 we show that for any $\Pb $ a solution for the Euler-Lagrange  equations, the operators $\Pi_\pm \Pb $ will also be a solution, regradless of whether the action is finite. In Theorem 5, we prove that if $\Pb $ is defined on the entire extended complex plane with finite action, then $\Pi_\pm \Pb $ will be as well. These two theorems are known in the vector case, \cite{Zakbook, GoldGrund2010} but new for the projector solutions. In Theorem 6, we prove the new results that given a finite action solution $\Pb $ for the Euler-Lagrange  equations,  using the covariant derivative defined for projectors, we can create two orthogonal subspaces of $\C$ which will be spanned by the set of rank-1 Hermitian projectors created from the raising and lowering operators. We prove that this set is mutually orthogonal and that the lowering operator will at some point annihilate the initial projector and hence there is some holomorphic projector associated with any solution to the $\CP$ model.

The new results of this section are similar to the analogous results, see \cite{ Zakbook} for the vector solutions but the proof using projectors has many advantages. By explicitly using the projective structure of $\CP,$ the equations are often more compact and the steps of the proof more transparent. Additionally, in the course of the proof, we have used the differential structures native to  $\CP$ and have proved several new results about such structures which can be used in contexts outside of the proof. An immediate example of this is the analysis of surfaces defined by the generalized Weierstrass formula for Immersions associated with  $\CP$ sigma models.

In Section 5, we consider surfaces immersed in $su(N)$ associated to $\CP$ model and apply the projective structures previously introduced to analyze the geometry of such surfaces.  In particular, we prove that any surfaces created from the Generalized Weierstrass formula for immersion using the conservation law derived from the Euler-Lagrange equations with finite action will be conformally parameterized by the Lagrangian density and will have finite surface area equal to the action of the model, Theorem 8.  Furthermore, we show that the Gaussian curvature, the norm of the mean curvature, the Willmore functional and the Euler-Poincar\' e characteristic can be written in terms of the Lagrangian and topological charge density which themselves can be written in terms of the norms of the covariant derivative of the projector.  We also show that the surface must be orthogonal to the sequence of rank-1 Hermitian projectors defined by repeated application of the raising and lowering operators. Next, we give necessary and sufficient conditions that a surface be associated with a $\CP$ model in terms of the characteristic polynomial and the differential. Finally, we illustrate the theoretical results by considering first arbitrary surfaces generated by holomorphic projectors and later an example based on the Veronese sequence and another with non-constant Gaussian curvature.

\section{The projector formalism}
To estabilsh our notation we start with a brief recalling of the mathematical grounds of the projector formalism. From computational purposes, we introduce a new sesquilinear product and compatible covariant derivative acting on the projector space. Next we describe certain characteristics of holomorphicity for projectors and basic properties of so called raising and lowering operators. Proofs have been put in the Appendix. 
\subsection{Hermitian, rank-1 projectors} 
The points of the complex coordinate space $\C$ are denoted by $\z=(z_1,.., z_n)^T$ and the Hermitian inner product on $\C,$ using $\dagger$ to denote Hermitian conjugation, as
\[\z^\dagger \w=(\z,\w)=\sum_{j=1}^N\overline{z}_jw_j, \qquad \z, \w \in \C.\]
The $(N-1)$ dimensional complex projective space $\CP$ is defined as a set of 1-dimensional subspaces of $\C;$ that is, all complex lines in $\C$ passing through a point $\z$ and the origin, together with the quotient topology.
Two points $\z,\w \in \C \diagdown \lbrace 0 \rbrace$ are in the same equivalence class if and only if there exists $a \ne 0 \in \mathbb{C}$ such that $\z=a\w.$ 

  By definition, a rank-1 Hermitian projector maps all of $\mathbb{C}^N$ on to some one dimensional subspace and so we can identify equivalence classes in $\CP$ with rank-1 Hermitian projectors. In  \ref{AA}, you can find the proof and explicit isomorphism between $\CP$ and the set of rank-1 projectors. 
 We say $\Pb $ is a rank-1 Hermitian projector if it has the following qualities, $ \Pb ^2=\Pb ,$ $ \Pb ^\dagger=\Pb $ and $tr(\Pb )=1$.

The most important properties of rank-1 Hermitian projectors in \\$GL_N(\mathbb{C}(\xi_+,\xi_-))$ are given below. 
\begin{property} From any rank-1 Hermetian projector $\Pb  \in GL_N(\mathbb{C}(\xi_+,\xi_-))$ we have the following identities
\begin{itemize}
 \item[1.1]  From any vector $\e,$ we can decompose $\Pb $ as proportional to the tensor product of $\Pb  \e$
\be \label{P11} \Pb  \e\otimes \e^\dagger \Pb =\e^\dagger \Pb  \e\Pb \ee
If $\Pb  \e$ is also assumed to be nonzero, we can solve for $\Pb $ and obtain,
\be \label{Pdecomp} \Pb =\frac{\Pb  \e\otimes \e^\dagger \Pb }{\e^\dagger \Pb  \e}.\ee 
\item[1.2] For any matrix $A,$ we have the following
\be \label{Ptrace}  \Pb  A\Pb =tr(\Pb  A)\ee
\item[1.3] We have the following differential constraint 
\be \label{dP} \partial_\mu \Pb =\partial_{\mu}\Pb  \Pb +\Pb \partial_{\mu}\Pb , \qquad \partial_\mu=\partial_{\xi_\mu}.\ee
\end{itemize}
\end{property}
The proofs can be found in  \ref{AB}.

\subsection{Sesquilinear Product and Covariant Derivative}
In this section, we introduce a sesquilinear product and covariant derivative, both depending on a fixed rank-1 Hermitian projector $\Pb $ and $\e \in \C$ a constant vector for which $\Pb  \e \ne 0.$ We do not chose $\e$ so that it is an eigenvector for $\Pb $ for reasons that will be discussed later. These structures are chosen in such a way that the covariant derivative $D_\mu ^{(\Pb )}$ is compatible with the sesquilinear product $\langle , \rangle _{(\Pb )}$ and at the same time the image is orthogonal to the argument.

\begin{definition} From a given $\Pb $ a rank-1 Hermitian projector and $\e$ a normalized vector with $\Pb  \e\ne 0,$ we define a sesquilinear product as 
$$\langle \ , \ \rangle: GL_N(\mathbb{C}(\xi_1, \xi_2))\times GL_N(\mathbb{C}(\xi_1, \xi_2))\rightarrow \mathbb{C}(\xi_1, \xi_2)$$
 by 
\be \label{<>} \langle \X,\Y\rangle _{(\Pb )}\equiv\frac{\e^\dagger \X \Y \e}{\e^\dagger \Pb  \e}.\ee
\end{definition}

Note that $\Pb  \e \ne 0$ which implies $\e^\dagger \Pb  \e=\e^\dagger \Pb  \e=(\Pb  \e)^\dagger \Pb  \e >0$ and so the inner product is well defined. Now, in general, $\e^\dagger \Pb  \e$ will be a scalar function, so we cannot normalize it to $1$ while retaining that $\e$ be an constant vector. 

This sesquilinear product is proportional to the Euclidean inner product on vectors induced by the transpose conjugate, i.e.  for $\x,\y \in \C,$ $\langle \x,\y\rangle_{Euc}=\x^\dagger \y.$ That is, if $\x\equiv \X \e$ and $\y\equiv \Y \e$ then 
$$\langle \X,\Y\rangle_{(\Pb )}=\frac{\e^\dagger \X^\dagger \Y \e}{\e^\dagger \Pb  \e}=\frac{\x^\dagger \y}{\e^\dagger \Pb  \e}.$$ 

We can use the  sesquilinear product to determine a semi-norm.
\be \label{norm} ||\X||_{(\Pb )}\equiv \sqrt{\langle \X,\X\rangle _{(\Pb )}}\ee
This is well defined since
\be \langle \X,\X\rangle _{(\Pb )}= \frac{(\X \e)^\dagger \X \e}{(\Pb  \e)^\dagger \Pb  \e} \geq 0.\ee
This norm is only semi-definite since any matrix orthogonal to $\e$ will have zero length, however we have the following property which is prove in \ref{AC}.
\begin{property} If $\X$ and $\Y$ are arbitrary rank one projectors and $\X e, \Y e \ne 0$, then 
\be\label{XY}  \langle \X,\Y\rangle _{(\Pb )}=0 \iff \X \Y=\0.\ee
\end{property}
 Although this is not a true norm, it will still satisfy many of the inequalities for norms. In particular, the Cauchy-Schwarz inequality and the triangle inequality which are proven in  \ref{AC}. 
From the remainder of the paper, we will omit the subscript  $\Pb $ in the sequilinear product unless there is ambiguity.

For computational convienence, in what follows we denote  $\mu$ an arbitrary direction in complex coordinates, $\xi_{\pm}=\xi_1\pm i\xi_2$  by $\mu=a \xi_++b \xi_-,$ $\overline{\mu}=a \xi_- +b \xi_+,$ for $ \ a,b \in \mathbb{R}.$ We identify the tangent space of $Gl_N(\mathbb{C}(\xi_1,\xi_2))$ with itself and denote $\vec{\mu}$  to be the vector in $Gl_n(\mathbb{C}(\xi_1,\xi_2)),$ generated by derivative $\frac{\partial}{\partial \mu}$ in the direction $\mu.$ The derivatives are represented as
\[ \partial_{+}\equiv \frac{\partial}{\partial \xi_{+}}=\frac{1}2\left(\frac{\partial}{\partial \xi_1}-i \frac{\partial}{\partial \xi_2}\right), \qquad \partial_{-}\equiv \frac{\partial}{\partial \xi_{-}}=\frac{1}2\left(\frac{\partial}{\partial \xi_1}+i \frac{\partial}{\partial \xi_2}\right).\]

 \begin{definition} We define a covariant derivative in the direction $\mu$ depending on the sesquilinear product given above, by 
\be \label{D} D_{\vec{\mu}}^{(\Pb )}=\frac{\partial}{\partial \mu}+\langle \Pb , \frac{\partial}{\partial \mu}\Pb \rangle_{(\Pb )}.\ee
\end{definition}

The covariant derivative in the basis directions $D_\pm ^{(\Pb )}$ are
\[D_{+}^{(\Pb )}\equiv D_{\xi_{+}}^{(\Pb )}=\p-\langle \Pb ,\p \Pb \rangle, \qquad D_{-}^{(\Pb )}\equiv D_{\xi_-}^{(\Pb )}=\bp -\langle \Pb ,\bp \Pb \rangle. \] 

Notice that this derivative depends on the direction $\vec{\mu},$ the operator $\Pb $ and implicitly the choice of $\e.$ We also note that the covariant derivatives $D_\pm ^{(\Pb )}$ in conjugate directions are not Hermitian conjugate to one another since 
$$(D_{+}^{(\Pb )})^\dagger-D_{-}^{(\Pb )}=\langle \Pb ,\partial_-\Pb \rangle-\langle \p \Pb , \Pb  \rangle$$
is not zero unless the norm $\e^\dagger \Pb  \e$ is a constant function.


\begin{property} 
 $D_{\vec{\mu}}^{(\Pb )}=\partial_ {\mu}+\langle \Pb , \partial_ {\mu}\Pb \rangle_{(\Pb )}$ is a covariant derivative. 
\end{property}
Proof. To prove this we need to show (see \cite{DoCarmo} for example), that for arbitrary matrices $\X, \Y,\vec{\mu}, \vec{\nu}  \in Gl_n(\mathbb{C}(\xi_1,\xi_2)), \ c,d \in \mathbb{R}$ and $\f: \mathbb{C} \rightarrow \mathbb{C},$
\bea &D^{(\Pb )}_{c\vec{\mu}+d\vec{\nu}}\X=cD^{(\Pb )}_{\vec{\nu}}\X+dD^{(\Pb )}_{\vec{\mu}}\X, &\text{linearity in direction} \\\
&D_{\vec{\mu}}^{(\Pb )}(\X+\Y)=D_{\vec{\mu}}^{(\Pb )} \X +D_{\vec{\mu}}^{(\Pb )} \Y,  &\text{linearity in argument}\\
&D_{\vec{\mu}}^{(\Pb )}(\f\X)=(\partial_{\vec{\mu}} \f)\X+\f D_{\vec{\mu}}^{(\Pb)} \X,  &\qquad \qquad  \text{Leibniz rule}
\eea 
These can be directly verified from the definition. 
From notational simplicity, we drop the superscript which denotes the dependence on $\Pb .$ This is particularly useful when we want to described repeated application of the covariant derivative, defined inductively by,  
$$D_\mu^k\X \equiv D_\mu (D_\mu^{k-1} \X), \qquad D_\mu^0\X\equiv \X .$$

The following are the most important properties of the covariant derivative. 
\begin{property} From the given rank-1 projector $\Pb $ and arbitary matrix functions $\X,\Y$ 
\begin{itemize}
\item[4.1]  $\Pb $ and  $D_\pm \Pb $ are orthogonal\\
\be \label{PDP} \langle \Pb ,D_{\pm}\Pb \rangle=0.\ee
\item [4.2] The covariant derivative, $D_\mu^{(\Pb)},$ is compatible with the sesquilinear product in the following sense
 \be \label{d<>}\frac{\partial}{\partial \mu} \langle \X,\Y\rangle =\langle D_{\overline{\vec{\mu}}}^{(\Pb )}\X,\Y\rangle+\langle \X,D_{\vec{\mu}}^{(\Pb )}\Y\rangle\ee
 \item[4.3] If $\X$ is a rank-1 Hermitian projector,  there exists a scalar function $\phi(\X)$ so that 
\be \label{DX} \partial_{\pm} \X \X \e=(D_{\pm}\X-\phi(\X) \X)\e\ee
 and in particular, if $\X=\Pb $ we have 
 \be \label{DP} \partial_{\pm} \Pb  \Pb  \e=D_{\pm}^{(\Pb )}\Pb  \e.\ee
 \item[4.4] The norm of $D_\pm \Pb $ is given by
 \be \label{normDP}|| D_+\Pb ||^2=tr(\p \Pb  \Pb  \bp \Pb ) \qquad ||D_-\Pb ||^2=tr(\bp \Pb  \Pb  \p \Pb )\ee
 \item[4.5]The commutator of covariant derivatives is given by
\be \label{Xpm} [D_{+},D_{-}]\X=(||D_+\Pb ||^2-||D_-\Pb ||^2) \X.\ee
 \end{itemize}
\end{property}
These properties are proven in \ref{AD}.

\subsection{(Anti-)Holomorphicity of vectors in the $\CP$  and their associated projectors}
Because of the phase invariance of $\CP,$ the definition of holomorphicity is differenct from other vector functions.  
\begin{definition}
We say a vector $\x$ in $\CP(\xi_+, \xi _-)$ is holomorphic if there is some non-zero scalar function $c(\xi_{+},\xi_{-}): \mathbb{R}^2\rightarrow \mathbb{C}$ for which 
$\bp (c\x)=\0.$
We call a rank-1 projector $\X$ a holomorphic projector if the vector which it projects onto is holomorphic.
\end{definition}
\begin{property} We have the following properties of holomorphic vectors and projectors.
\begin{itemize}
\item[5.1] A vector $x$ in  $\CP$ is holomorphic if and only if $\bp \x$ is proportional to $ \x.$
\item[5.2] A Hermitian, rank-1 projector $\X$ is holomorphic if and only if 
\be \label{projholo} \X\partial_+ \X=\partial_-\X\X=\0. \ee 
\item[5.3] If $\X=\Pb $ the projector which was used to define the covariant derivative \eqref{D} then $\Pb $ is holomorphic if and only if 
\be D_-\Pb  \e=0.\ee
\end{itemize}
\end{property}
The proofs are in \ref{AE}.\\
The analogous statements for anti-holomorphicity follow directly. That is, $\X$ is antiholomorphic if the projector which it projects onto is antiholomorphic and this is equivalent to 
\be \X\partial_- \X=\partial_+\X\X=\0.\ee
If $\X=\Pb $ then it is antiholomorphic if and only if $D_+\Pb  \e=0$

\subsection{Basic Properties of Raising and Lowering Operators}
We have the following raising and lowering operators from previous literature \cite{GoldGrund2009, GoldGrund2010} 
\be \Pi_{\pm}(\X)\equiv \Bigg{\lbrace} \begin{array}{cc} \frac{\partial_\pm \X \X\partial_\mp \X}{tr(\partial_\pm \X \X\partial_\mp \X)} & \partial_\pm \X \X\partial_\mp \X \ne \0\\
\0 &\partial_\pm \X \X\partial_\mp \X = \0\end{array}, \ee
These raising and lowering operators are analogous to the raising and lowering operators for the vector model, see for example \cite{Zakbook}.
These operators are defined in such a way that the $\Pi_{\pm}(\X)$ is a rank-1 Hermitian projector whenever $\X$ is a rank-1 Hermitian projector and so that $\X$ is orthogonal to $\Pi_{\pm}\X$ and hence the operation is similar to a Gram-Schmidt orthogonalization procedure. We prove these results in the following properties. 
\begin{property} From $\X$ a Hermitian rank-1 projector we have 
\begin{itemize}
 \item[6.1]$\Pi_{\pm}\X$ is a Hermitian projector which projects onto the vector $ \partial_\pm \X \X e.$
 Thus, either $\Pi_{\pm}\X$ is $\0$ or it is of rank-1 and can be written as 
\be \label{PX} \Pi_{\pm}\X=\frac{ \partial_\pm \X \X \e\otimes \e^\dagger \X \partial_\mp \X}{\e^\dagger \X\partial_\mp \X\partial_\pm \X \X \e}.\ee

\item[6.2]  $\X$ is orthogonal to $\Pi_{\pm}\X.$
\item[6.3] The raising and lowering operators acting on $\Pb $ can be written in terms of the covariant derivative 
\be\label{rlcov} \Pi_{\pm}\Pb =\frac{D_{\pm}\Pb  \e\otimes \e^\dagger (D_{\pm}\Pb )^\dagger}{(D_{\pm}\Pb  \e)^\dagger D_{\pm}\Pb  \e}\ee
\item[6.4]  $\X$ is holomorphic if and only if 
\be \label{PiXanal}\Pi_{-}\X=\0.\ee
\end{itemize}
\end{property}
The proofs can be found in \ref{AF}. Again, there is an analogous Property 5.4 for anti-holomorphicity. That is, $\X$ a rank one projector is antiholomorphic if and only if $\Pi_+\X=\0.$ Finally, we can inductively define repeated application of the raising operator by 
$$\Pi^k_\pm(\X)=\Pi_\pm(\Pi_\pm^{k-1}(\X)), \qquad \Pi_\pm^0(\X)=\X,\ k \in \mathbb{N}.$$

\section{The Action and the  Euler-Lagrange Equation}
In this section we introduce the basic quanties and equations of the $\CP$ model and show that they can be expressed in terms of the gauge invariant covariant derivative \eqref{D}. From the $\CP$  model, we can write the \textbf{ Lagrangian density} in terms of projectors in two equivalent ways. 
\be \label{lagrangian1} L(\Pb )=  tr(\partial_+ \Pb \bp \Pb ) \ee

This is the form most often found in the literature though without a factor of $2$ as in some of the references \cite{ZakMik1979, Mik1986, Zakbook}.  The following form of the Lagrangian looks identical to the one for the vector model but the definition of the covariant derivative and inner product are different. 

\be \label{lagrangian2} L(\Pb )=
||D_{+}^{(\Pb )}\Pb ||_{(\Pb )}^2+||D_{-}^{(\Pb )}\Pb ||_{(\Pb )}^2 \ee

Here we write the dependence of the inner product on $\Pb $ explicitly as a reminder that the covariant derivatives and the semi-norm are both dependent on $\Pb .$
The second formulation of the Lagrangian density is a direct consequence of the action as written in terms of the vector $f$ onto which $\Pb $ projects, ie. $\f=\Pb  \e.$ 

\begin{theorem}\label{laeq}
The two formulations of the Lagrangian density are equivalent and in particular, the second formulation does not depend on the choice of $\e.$
\end{theorem}
Proof. This follows immediate from \eqref{normDP} which gives the semi-norm of $D_\pm \Pb $ in terms of trace. 
\QED

There is another quantity associated with the model, called \textbf{the topological charge density}  and it is given by
\be \label{chargedensity} q(\Pb )=||D_+\Pb ||^2-||D_-\Pb ||^2=tr(\Pb \p \Pb  \bp \Pb -\Pb \bp \Pb \p \Pb )\ee
It can be shown that this quantity is a total divergence and hence the integral 
\be\label{charge} Q(\Pb ) =\frac 1 \pi \int_{\mathbb{R}^2}q(\Pb ) d\xi_1 d\xi_2\ee
is an integer \cite{DinZak1980prop}. 
The topological charge density has already appeared in this paper in  \eqref{Xpm} the commutator of covariant derivatives,
\[ [D_{+},D_{-}]\X =q(\Pb ) \X, \qquad \forall \X \in Gl_N(\mathbb{C})(\xi_+, \xi_-)\]

 The \textbf{ action} for the $\CP$ sigma model on the Euclidean plane $\mathbb{R}^2$ is given by 
\be \label{action} S(\Pb )=\int_{\mathbb{R}^2}\mathcal{L}(\Pb )d\xi_1d\xi_2.\ee

There are also several equivalent formulations of the \textbf{ Euler-Lagrange equations} for the fixed points of the action.  A rank-1 projector $\Pb $ is a fixed point of the action if and only if one of the following hold:

\be \label{EL} [\partial_+ \bp \Pb ,\Pb ]=\0\ee
which is equivalent to 
\be \label{el2} \partial_+ [\bp \Pb ,\Pb ] +\bp[\partial_+ \Pb ,\Pb ]=\0.\ee

The other formulation of the Euler-Lagrange  equations is via
\be \label{el3} \left(D_{+}D_{-}(\Pb )+||D_{-}\Pb ||^2\Pb \right)\e=\0\ee
which is equivalent via \eqref{Xpm} to 
\be \left(D_{-}D_{+}(\Pb )+||D_{+}\Pb ||^2\Pb \right)\e=\0.\ee

\begin{theorem}
Euler -Lagrange Equations \eqref{EL} and \eqref{el3} are equivalent.
\end{theorem}
Proof, We shall prove that Euler-Lagrange formulations \eqref{EL} and \eqref{el3} both are equivalent to the statement $\p \bp \Pb  \Pb $ is proportional to $\Pb .$  By the projective property of $\Pb $, $\p \bp \Pb  \Pb $ is proportional to $\Pb $ implies that the constant of proportionality be  $tr(\Pb \p\bp \Pb )$ and so   
\be\label{elm} \p \bp \Pb  \Pb =tr(\Pb \p\bp \Pb ) \Pb . \ee
Assume $\Pb $ satisfies \eqref{EL}, then  $(1-\Pb ) \p \bp \Pb  \Pb =\0$
which imples that $\Pb $ satisfies \eqref{elm}.
On the other hand,  assume  $\Pb $ satisfies \eqref{elm}. We take the  Hermitian conjugation of \eqref{elm} to obtain
\[  \Pb  \p \bp \Pb  =tr(\Pb  \p \bp \Pb ) \Pb \] 
which, together with \eqref{elm}, imples
\[  [\p \bp \Pb , \Pb ]=\0 \]
So, \eqref{EL} is equivalent to \eqref{elm}. To prove the equivalence of \eqref{el3} and \eqref{elm}, we can use \eqref{DP} to rewrite  \eqref{el3} as
\[ \p \bp \Pb  \Pb  \e +\bp \Pb  \p \Pb  \e -\frac{\e^\dagger \Pb  \p \Pb  \e}{\e^\dagger \Pb  \e}\bp \Pb  \Pb  \e +tr(\bp \Pb \p \Pb )\Pb  \e=0. \]
Using \eqref{Ptrace} and \eqref{dP}, we obtain
\bea\bp \Pb  \p \Pb  \e&=&\Pb  \bp \Pb  \p \Pb  \Pb  \e+\bp \Pb  \Pb  \p \Pb  \e\nn
&=&tr(\Pb  \bp \Pb  \p \Pb  )\Pb  \e+\frac{\e^\dagger \Pb  \p \Pb  \e}{ \e^\dagger \Pb  \e}\bp \Pb  \e\eea
which show that Euler-Lagrange  \eqref{el3} is equivalent to 
$$\p \bp \Pb  \Pb  \e + tr(\Pb  \bp \Pb  \p \Pb )\Pb  \e +tr(\Pb  \p \Pb  \bp \Pb )\Pb  \e=0$$
and so $\p \bp \Pb  \Pb  \e=b \Pb \e$ for some scalar function $b(\xi_+, \xi_-).$ As before, if we multiply both sides of the equation on the left by $\Pb $, we see that $b=tr(\Pb  \p\bp \Pb ).$ 

Hence, 
$$\left(D_{+}D_{-}(\Pb )+||D_{-}\Pb ||^2\Pb \right)\e=0 \iff \p \bp \Pb  \Pb  \e =tr(\Pb  \p \bp \Pb )\Pb \e.$$ 
But also, for any vector $\v,$ there exists some scalar function $c=c(\xi_+, \xi_-)$ such that $\Pb  \v=c\Pb  \e$ and so 
$$ (\p \bp \Pb ) \Pb  \v = c(\p \bp \Pb ) \Pb  \e= c tr(\Pb  \p \bp \Pb )\Pb  \e=tr(\Pb  \p \bp \Pb )\Pb  \v.$$
Thus, we see that Euler-Lagrange equation  \eqref{el3} is equivalent to \eqref{elm}, that is
$$\left(D_{+}D_{-}(\Pb )+||D_{-}\Pb ||^2\Pb \right)\e=0 \iff \p \bp \Pb  \Pb  =tr(\Pb  \p \bp \Pb )\Pb .$$

Hence, we have shown that the Euler-Lagrange equations \eqref{EL} and \eqref{el3} are equivalent to the requirement that $\p \bp \Pb  \Pb  $ is proportional to $\Pb , $ and in particular are equivalent to each other.  \QED

From $\Pb $ a solution to the Euler-Lagrange equations, the rasing and lowering operators are mutual inverses. This theorem was proven in \cite{GoldGrund2009, GoldGrund2010} and we give the proof in \ref{AG}
\begin{theorem} If $\Pb $ satisfies the Euler-Lagrange equation \eqref{EL}, then 
\[\Pi_\pm(\Pi_\mp \Pb )=\Pb ,\qquad \Pi_\mp(\Pi_\pm \Pb  )=\Pb \] 
whenever  $\Pi_\pm \Pb\ne \0$.\end{theorem}  

\section{Projector solutions for the equations of motion}
In this section we prove some results about projectors which satsify the equations of motion \eqref{EL}. Namely, we can use the raising and lowering operators to build a sequence of projectors which will be solutions to the Euler-Lagrange equation and furthermore all solutions are in one such family generated by a holomorphic projector. 

The following theorem has been proven for the vector case with the raising and lowering operators associated with them, see \cite{Zakbook} and references therein. We prove it here directly for projectors. 
\begin{theorem}\label{Pi e-l}
If $\Pb $ is rank-1 Hermitian projector which is a solution of the Euler-Lagrange equations \eqref{EL} then $\Pi_{\pm}\Pb $ will be as well. 
\end{theorem}
Proof. We shall prove that $\Pi_+\Pb $ satisfies Euler-Lagrange \eqref{EL}, the corresponding proof for $\Pi_-\Pb $ follows by direct analogy. Assume that $\Pb $ is a rank-1 projector satisfying \eqref{EL}. We shall now prove that there exists some complex function $b$ so that 
\be \p \bp \Pi_+\Pb  \Pi_+\Pb =b \Pi_+\Pb \ee
which, as disscussed above, is equivalent to Euler-Lagrange  equation, for example \eqref{EL}. 

The computation is purely algebraic and uses the properties of rank-1 projectors and the fact that $\Pb $ is a solution of \eqref{EL}.
We  uses the identity 
$$tr(\Pb \p \bp \Pb )=-tr(\p \Pb \Pb \bp \Pb +\bp \Pb  \Pb  \p \Pb )$$
which is a differential consquence of $tr(\Pb )=1$ and $\Pb ^2=\Pb .$ Also, \eqref{EL} implies
$$tr(\p \bp \Pb  \partial_\pm \Pb )=tr(\p \bp \Pb  \Pb  \partial_\pm \Pb ) =0 .$$
These identies are used to compute
\be \p \bp \Pi_+\Pb \Pi_+\Pb  =b\Pi_+\Pb \ee

where, 
\[ b\equiv \frac{tr(\Pb \p \Pb  \p \Pb )tr(\Pb  \bp \Pb  \bp \Pb ) +tr(\Pb  \p \Pb  \bp ^2\Pb  \p \Pb )+tr(\Pb  \bp^2\p \Pb  \p \Pb )}{tr(\p \Pb  \Pb  \bp \Pb )}\]
\[+\frac{\p tr(\p \Pb  \Pb  \bp \Pb )\bp tr(\p \Pb  \Pb  \bp \Pb )-\p\bp tr(\p \Pb  \Pb  \bp \Pb )tr(\p \Pb  \Pb  \bp \Pb )}{tr(\p \Pb  \Pb  \bp \Pb )^2}\]
\be-tr(\p \Pb  \Pb  \bp \Pb ).\ee

Thus, $\Pi_+\Pb $ is a solution to the equations of motion. Finally, we note that by  \eqref{PX}, $\Pi_+\Pb $ is a rank-1 Hermitian projector so the theorem is proven. 
\QED

Let $\overline{\mathbb{C}}\equiv \mathbb{C}\cup \lbrace \infty \rbrace $ denote the extended complex plane. We shall now prove that if $\Pb $ is a finite action solution defined on $\overline{\mathbb{C}},$ then $\Pi_\pm \Pb $ will be as well.  
\begin{theorem}\label{Pi action}
 If $\Pb $ is solution to the Euler-Lagrange equation \eqref{EL} in $C^{M+2}(\overline{\mathbb{C}})$ for $ 0\leq M\in \mathbb{N}$ and has finite action \eqref{action}, then the the operator $\Pi_{\pm}\Pb $ will also be  a solution to the Euler-Lagrange equation in $C^{M+1}(\overline{\mathbb{C}})$ with finite action.
\end{theorem}
Before we prove this theorem, we will prove a lemma about the inversion of the projector about the unit sphere. 
\begin{definition} From any matrix function $\X(\xi_+, \xi_-) ,$ we define the inversion about the unit sphere as 
\be \label{hat}\hat{ \X} (w_+,w_-) \equiv \X(\xi_{+},\xi_-), \qquad \xi_{+} \equiv w_{+}^{-1}, \ \xi_{-} \equiv w_{-}^{-1}.\ee
\end{definition}

\begin{lemma}\label{hatP}
If $\Pb (\xi_+, \xi_-)$ is solution to the Euler-Lagrange equations \eqref{EL} in $C^{M+2}(\overline{\mathbb{C}})$ for $ 0\leq M\in \mathbb{N}$ and with finite action \eqref{action} then $\hat{ \Pb } (w_+,w_-)$ will also be a solution to the Euler-Lagrange equations in $C^{M+2}(\overline{\mathbb{C}})$ with finite action.
\end{lemma}
Proof. Assume $\Pb (\xi_+, \xi_-)$ is solution of the Euler-Lagrange equations in $C^{M+2}(\overline{\mathbb{C}})$ for $ 0\leq M\in \mathbb{N}$ and has finite action. Now, we define a new projector $\hat{\Pb }$ by inverting $(\xi_+, \xi_-)$ across the unit circle via \eqref{hat}. Because of the change of variables, we must return to the more cumbersome notation to make explicit the set of variables to which we are referring. From notational convience we write $\Pb =\Pb (\xi_+, \xi_-), \ \hat{\Pb }=\hat{\Pb }(w+,w_-).$
We claim that $\hat{\Pb }$ will also be a rank-1 projector solution in $C^{M+2}(\overline{\mathbb{C}})$ with finite action. To show this, we note that 
\be \label{COV} \partial_{w_{\pm}}\hat{\Pb }=-\frac 1{w_{\pm}^2}\partial_{\xi_{\pm}}\Pb , \qquad D_{w_{\pm}}^{(\hat{\Pb })}{\hat{\Pb }}=-\frac{1}{w_{\pm}^2}D_{\xi_{\pm}}^{(\Pb )} \Pb ,\ee 
so that any of the forms of the Euler-Lagrange equation for $\hat{\Pb }$ will be the Euler-Lagrange equation for $\Pb $ multiplied by $w_+^{-2}w_-^{-2}.$ From example, 
$$[\partial_{w_+}\partial_{w_-}\hat{\Pb },\hat{\Pb }]=\frac{1}{w_+ ^2w_-^2}[\partial_{\xi_+}\partial_{\xi_-}\Pb ,\Pb ]=0.$$
Similarly, 
$$S(\Pb )=\int_{\mathbb{R}^2} <D_{w_{+}}^{(\hat{\Pb })}{\hat{\Pb }},D_{w_{+}}^{(\hat{\Pb })}{\hat{\Pb }}>_{(\hat{\Pb })}+<D_{w_{-}}^{(\hat{\Pb })}{\hat{\Pb }},D_{w_{-}}^{(\hat{\Pb })}{\hat{\Pb }}>_ {(\hat{\Pb })}dw_+dw_-$$
$$=\int_{\mathbb{R}^2}<D_{\xi_{+}}^{(\Pb )}\Pb ,D_{\xi_{+}}^{(\Pb )}\Pb >_{(\Pb )}+<D_{\xi_{-}}^{(\Pb )}\Pb ,D_{\xi_{-}}^{(\Pb )}\Pb >_{(\Pb )}d\xi_+d\xi_-=S(\hat{\Pb }).$$
So, if $\Pb $ is a finite action solution on the extended complex plane $\overline{\mathbb{C}}$, $\hat{\Pb }$ will be as well. \QED

Proof of Theorem \ref{Pi action}. Because of Theorem \ref{Pi e-l}, $\Pi_{\pm}\Pb $ will satisfy the equations of motion and so it suffices to prove that $\Pi_{\pm}\Pb $ has finite action. That is
$$S(\Pi_{\pm}\Pb )=\int_{\mathbb{R}^2}tr(\partial_+ \Pi_{\pm}\Pb \bp \Pi_{\pm}\Pb )dx_{+}dx_{-}<\infty$$

To prove this, we shall transform the integral into two integrals of real analytic functions over the unit disk and so it will be finite. In order to do this we use a result of A. Friedman, that given any real set of nonlinear elliptic partial differential equations which is analytic in the arguments, the solutions which are at least twice differentiable must be real analytic \cite{Friedman}. In our particular case, we see that the real and imaginary parts of the  components of our projector $\Pb $ , $\hat{\Pb }$ can be decomposed into real functions and then the sigma model can be described as a set of elliptical partial differential equations for these components, see for example H. J. Borchers and W. D. Garber \cite{BorchersGarber}. Thus, we can apply the results of Friedman to see that if such functions are at least twice integrable, they will be real analytic. 

Now, since  the components of $\Pb $ have real and imaginary parts which are analytic, we can use \eqref{PX} to show that the components of $\Pi_{\pm}\Pb $ also have real and imaginary parts which are analytic. Recall
$$\Pi_{\pm}\Pb =\frac{ \partial_{\pm}\Pb  \Pb  \e\otimes \e^\dagger \Pb  \partial_{\mp} \Pb }{\e^\dagger \Pb \partial_{\mp}\Pb \partial_{\pm} \Pb \Pb  \e},$$
and so $\Pi_{\pm}\Pb $ will have real analytic components as long as  the denominator of $\Pi_{\pm}\Pb $ given by 
\be \label{denom}  \e^\dagger \X\partial_{\mp}\X\partial_{\pm} \X \X \e=(\X\partial_{\pm} \X \X \e)^\dagger \X\partial_{\pm} \X \X \e,\ee
is not zero.  However, we see that this term \eqref{denom}  is $0$ if and only if $\Pi_{\pm}\Pb =0$ in which case the assertion is trivial. Thus, the real and imaginary parts of the components of $\Pi_{\pm}\Pb $ are real analytic and hence the Lagrangian density will be as well. 

We will now split the area of the integral into the interior and exterior of the unit disk. Since our integrand is analytic on the whole plane, including the unit circle, the boundary of our domains will not affect the integral. By the analyticity of $\mathcal{L}(\Pi_\pm \Pb ), $ the integral of the Lagrangian density over the unit disk becomes an integral of real analytic functions over a compact domain so is finite,  
\be \label{intP1}\int_{|\xi_+|\leq 1}\mathcal{L}(\Pi_{\pm}\Pb )d\xi_1\xi_2= \int_{|\xi_+|\leq 1}tr(\partial_+ \Pi_{\pm}\Pb \bp \Pi_{\pm}\Pb ) d\xi_1\xi_2 <\infty.\ee 
Next, we integrate $\mathcal{L}(\Pi_{\pm}\Pb )$ over the exterior of the unit disk $|\xi_+|\geq 1$ by making a change of variables $\xi_\pm=w_\pm^{-1}$  and integrating over the interior of the unit disk $|w_+|\leq 1.$

To do this, we make use of Lemma \ref{hatP} which says if $\Pb (\xi_{+},\xi_-)$ is a rank-1 projector solution which is at least $C^2(\overline{\mathbb{C}})$ and with finite action then the function obtained by inverting the coordinates across the unit circle will also be a solution to the sigma model with finite action and, in particular, the components of $\hat{\Pb }$  have real and imaginary parts which are analytic.

Next, we invert $\Pi_{\pm}\Pb $ about the unit circle using \eqref{COV} as
\be \widehat{\Pi_{\pm}\Pb } =\frac{\partial_{w_{\pm}}\hat{\Pb }\hat{\Pb }\partial_{w_{\mp}}\hat{\Pb }}{tr(\partial_{w_{\pm}}\hat{\Pb }\hat{\Pb }\partial_{w_{\mp}}\hat{\Pb })},\ee
whose components  have real and imaginary parts which are analytic $w_+,w_-.$
Finally, we use the change of variables again to compute
$$\partial_{\xi_{\pm}}\Pi_{\pm}\Pb =-w_{\pm}^2\partial_{w_{\pm}}\widehat{\Pi_{\pm}\Pb }$$
and so we have an integral of real analytic functions over a compact domain,
\bea \label{extP1}\int_{|\xi_+|\geq 1}\mathcal{L}(\Pi_{\pm}\Pb )d\xi_1\xi_2= \int_{|\xi_+|\geq 1}tr(\partial_{\xi_+} \Pi_{\pm}\Pb \partial_{\xi_-} \Pi_{\pm}\Pb )d\xi_+d\xi_-\nn
=\int_{|w_+|\leq 1}tr(\partial_{w_+} \widehat{\Pi_\pm \Pb }\partial_{w_-}\widehat{\Pi_{\pm}\Pb })dw_+dw_-<\infty\eea
Thus, have shown that the integral of $\mathcal{L}(\Pi_{\pm}\Pb )$ over both the interior  \eqref{intP1} and the exterior \eqref{extP1} of the unit disk is finite and so the integral over the entire plane is finite as well. That is, 
$$S(\Pi_{\pm}\Pb )=\int_{\mathbb{R}^2}tr(\partial_{w_+} \widehat{\Pi_\pm \Pb }\partial_{w_-}\widehat{\Pi_{\pm}\Pb })dw_+dw_-<\infty.$$ \QED

Now that we know that we can use the raising and lowering operators to build sequences of projector solutions, we show that the sequence of projectors must terminate, that is there is a highest and lowest projector which is annihilated by the raising and lowering operators, respectively. Such projectors are anti-holomorophic or  holomorphic respectively and so we can see that any projector solution to the Euler-Lagrange equations with finite action belongs to some family of such projectors generated by a holomorphic or anti-holomorphic projector. We prove these results in the following theorem.

\begin{theorem}\label{bigone}
From any $\Pb :\overline{\mathbb{C}} \rightarrow Gl_n(\mathbb{C})$ a finite action, \eqref{action}, solution of the Euler-Lagrange equations, \eqref{EL},  which is a rank-1 Hermitian projector whose components are differentiable at least $N+2$ times, there exists a holomorphic projector $\Pb _0$  and a anti-holomorphic projector $\Pb _{m+n}$ such that \Pb  can be obtained by successive application of the raising operator to $\Pb _0$ or lowering operator to $\Pb _{m+n},$ i.e. there exist $n, m$ with $0\leq n+m \leq N-1$ such that
\[\Pb =\Pi_{+}^n\Pb _0, \qquad \Pb =\Pi_{-}^m \Pb _{m+n}, \qquad  \Pb _0=\Pi_{-}^{m+n}\Pb _{m+n}.\]
\end{theorem}

We shall prove this in the following lemma's.  
\begin{itemize}
\item[  \ref{1}] The vector spaces $H,H^\perp \subset \C$ as given by 
$$H\equiv \lbrace (D_{-}^k\Pb )\e | k\ge 0\rbrace, \qquad H^\perp\equiv \lbrace (D_{+}^k\Pb )\e |k\ge 1\rbrace$$ are in fact mutually orthogonal with respect to the sesquilinear product \eqref{<>}. 
\item[  \ref{2}] Let $n+1$ be the dimension of $H$ and $m$ the dimension of $H^\perp$. We shall show that the projectors $\Pi_{-}^k\Pb , 0\leq k$ project in to $H$ and that the projectors $\Pi_{+}^k\Pb , \ 1\leq k$ project in to $H^\perp.$ Furthermore the images of $\Pi_{-}^k\Pb , 0\leq k\leq n$ form a basis for $H$ and $\Pi_{+}^k\Pb , 1\leq k\leq m$ form a basis for $H^\perp.$
\item[ \ref{3}] The projectors $\Pi_{\pm}^k\Pb , k\in \mathbb{N}$ are mutually orthogonal and we have a decomposition of the operator 
$$I_{n+m+1}=\Pi_{-}^{n}\Pb +...+\Pi_{-}\Pb +\Pb +\Pi_{+}\Pb +...+ \Pi_{+}^{m}\Pb $$ 
which acts as the identity on $H\cup H^\perp.$
 \end{itemize}

From each step, we assume that $\Pb $ is a rank-1 Hermitian projector as a function of two variables on the Riemann two sphere which is a solution to the Euler-Lagrange equations \eqref{EL} with finite action \eqref{action}. Without loss of generality, we chose $\e$ so that $\Pi_\pm^k\Pb  \e$ is not identically $0$ unless the operator $\Pi_\pm^k\Pb =\0$ and  $\Pb  \e\ne 0$ for any $(\xi_+, \xi_-).$ Notice, that because $\Pb $ is orthogonal to $\Pi_\pm(\Pb ),$ we cannot chose the vector $\e$ as an eigenvector for $\Pb , $ since then it would be annihilated by the operators $\Pi_\pm^k\Pb .$

\begin{lemma}\label{1}
The vector subspaces $H, H^\perp\subset \C$ defined by
$$H\equiv \lbrace (D_{-}^k\Pb ) \e |0\leq k  \rbrace, \qquad H^\perp\equiv \lbrace (D_{+}^k\Pb ) \e|1\leq k \rbrace$$ are in fact mutually orthogonal, in the Euclidean norm.
\end{lemma}
Proof. We prove this by considering the functions 
$$A_{ij}\equiv\langle D_{-}^i\Pb ,D_{+}^j\Pb \rangle_{(P)}=\e^\dagger P\e \langle D_{-}^j\Pb  \e, D_{-}^k\Pb \e\rangle_{Euc}, \qquad i,j\in \mathbb{N}$$
and showing that they are identically $0$ for all $i+j>0$ and $i\leq N, j \leq N.$ Recall, $\e^\dagger P\e$ was assumed to be not identically $0..$ We note that the result depends on the fact that 
$D_{+}^j\Pb $ must be at least twice differentiable. In particular, if $\Pb $ is assumed to be at least $N+2$ differentiable, then $D_{+}^j\Pb $ will be at least twice differentiable for all $0\leq j \leq N.$

We prove this by induction on the quantity $n=i+j.$ From $n=1,$ by  \eqref{PDP} we see
$$0=A_{0,1}=\langle \Pb ,D_{+}\Pb \rangle,\qquad 0=A_{1,0}=\langle D_{-}\Pb ,\Pb \rangle.$$
Next, assume $A_{i,j}=0$ for all $0<i+j<n.$ We shall prove that $A_{i,j}=0$ for any $i+j=n.$ This can be even further simplified since for $i+j=n$
\bea \label{indexshift} A_{i,j}&=&\langle D_{-}^i\Pb ,D_{+}^{j}\Pb \rangle=\partial \langle D_{-}^{i-1}\Pb ,D_{+}^{j}\Pb \rangle-\langle D_{-}^{i-1}\Pb ,D_{+}^{j+1}\rangle \nn
&=&\partial A_{i-1,j}-A_{i-1,j+1}=-A_{i-1,j+1}\eea
here we have used the compatibility of the covariant derivative with the sesquilinear product,  \eqref{d<>}, and  the induction assumption that $A_{i',j'}=0$ for $i'+j'<n.$ By continued application of this procedure, we can see that showing 
$A_{0,n}=0$ implies that $A_{i,j}=0$ for all $i+j=n>0.$ We note here, that we could instead have continued the proof with $A_{n,0}$ and the proof would have been  identical under an interchange of $\pm.$

Now, we consider $A_{0,n}=\langle \Pb ,D_{+}^n\Pb \rangle$ and show that it is identically equal to $0$ for $n,$ under the induction assumption that $A_{i,j}$ is identically equal to $0$ for all $0<i+j<n.$ We shall prove this first by showing that the function $A_{0,n}$ is analytic and then use the finiteness of the action to prove that it is bounded and therefore identically equal to $0$.  

Consider, 
\[ \bp A_{0,n}=\langle D_{+}\Pb ,D_{+}^{n}\Pb \rangle+\langle \Pb , D_{-}D_{+}^n\Pb \rangle \]
We shall prove that this is identically $0$ by defining the scalar function
\be \label{B} B\equiv \langle D_{+}\Pb , D_{+}\Pb \rangle\ee 
and showing  that
\[ \langle D_{+}\Pb ,D_{+}^{n}\Pb \rangle=\partial_+^{n-1}B,\qquad \langle \Pb ,D_{-}D_{+}^{n}\Pb \rangle=-\partial_+^{n-1}B, \ \forall n\geq 1.\] 
where we define $\partial_\pm ^0B\equiv B$. 
From the first equality $,\langle D_{+}\Pb ,D_{+}^{n}\Pb \rangle=\partial_+^{n-1}B,$ we prove it inductively, for $1\leq k\leq n$ with the base step being true by definition. 
Next, assume that for $k-1 \geq 1,$ $\langle D_{+}\Pb ,D_{+}^{k-1}\Pb \rangle=\partial_+^{k-2}B.$ We shall prove the identity for $2 \leq k \leq n.$ That is, we see that
\be \langle D_{+}\Pb ,D_{+}^{k}\Pb \rangle=\partial_+ \langle D_{+}\Pb ,D_{+}^{k-1}\Pb \rangle-\langle D_{-}D_{+}\Pb ,D_{+}^{k-1}\Pb \rangle =\partial_+^{k-1}B.\ee
Here we used  \eqref{d<>} for the first equality and the induction assumptions $\langle D_{+}\Pb ,D_{+}^{k-1}\Pb \rangle=\partial_+^{k-2}B,$ to show that $\partial_+ \langle D_{+}\Pb ,D_{+}^{k-1}\Pb \rangle =\partial_+^{k-1}B$. Also we used $A_{0,n'}=0$ for $n'<n$ and the fact that $\Pb $ satisfies Euler-Lagrange \eqref{el3} to show that $\langle D_{-}D_{+}\Pb ,D_{+}^{k-1}\Pb \rangle=0.$
Thus, we have proven that for all $1\leq k\leq n, \ \langle D_{+}\Pb ,D_{+}^{k}\Pb \rangle=\partial_+^{k-1}B.$  

Next, we shall prove, that 
\be \label{ind3} \langle \Pb ,D_{-}D_{+}^{n}\Pb \rangle=-\partial_+^{n-1}B.\ee
To do this, we prove inductively, that 
\be \label{ind2}\langle \Pb ,D_{-}D_{+}^{k}\Pb \rangle =\langle \Pb ,D_{+}^{k-1}D_{-}D_{+}\Pb \rangle, \qquad  1\leq k\leq n \ee
This is trivial for $k=1.$ Now, assume its true for $k-1$ that is assume 
\be \label{ind2 k-1} \langle \Pb ,D_{-}D_{+}^{k-1}\Pb \rangle =\langle \Pb ,D_{+}^{k-2}D_{-}D_{+}\Pb \rangle\ee 
We use  \eqref{Xpm} which says $[D_+,D_-]\X=q\X$ for $q$ the topological charge density. 
Using this property and the induction assumptions \eqref{ind2 k-1} and  $A_{i,j}=0$ for $i+j<n,$ we obtain
\bea \langle \Pb ,D_{-}D_{+}^{k}\Pb \rangle &=&\langle \Pb ,(D_{+}D_{-}-q)D_{+}^{k-1}\Pb \rangle =\langle \Pb ,D_{+}D_{-}D_{+}^{k-1}\Pb \rangle -q A_{0,k-1}\nn
&=&\langle \Pb ,D_{+}^{k-1}D_{-}D_{+}\Pb \rangle \eea
and so we have proven \eqref{ind2} by induction.
 
Next, we use the fact that $\Pb $ is a solution to the Euler-Lagrange equations,  the Leibniz rule for covariant derivatives and the induction assumption that $A_{i,j}=0$ for $i+j <n$ to show 
\bea \langle \Pb ,D_{+}^{n-1}D_{-}D_{+}\Pb \rangle &=&=\langle \Pb ,D_{+}^{n-1}(B \Pb )\rangle =\langle \Pb ,(\partial_+^{n-1}B)\Pb \rangle =\partial_+^{n-1}B .\eea
Thus, we see that \eqref{ind3} holds and so we have proven that 
$$ \bp A_{0,n}=0.$$
Next, we use  \eqref{indexshift} and the Cauchy-Schwarz inequality for the semi-norm to show for $i+j=n$
\be |A_{0,n}|=|A_{i,j}|=|<D_{-}^i\Pb ,D_{+}^j\Pb >|\leq ||D_{-}^i\Pb ||_{(\Pb )}||D_{+}^j\Pb ||_{(\Pb )}.\ee
It suffices to show that the function $||D_{\pm}^j\Pb||_{(\Pb )} \in L^2(\mathbb{R}^2)$ for any $j\geq 1$ since by H\"older's inequality the product of two $L^2(\mathbb{R}^2)$ functions is $L^1(\mathbb{R}^2).$
Again, we recall that we are proving this by induction on $n$ with the case $n=1$ already shown and so we only consider  $n\geq 2$ and hence we can chose  $i,j \geq 1.$

In order to compute the integral we use an argument similar to the one in the proof of Theorem \ref{Pi action} where we split the integral into two part and write each part as an integral for real analytic functions over unit disks. 

 We know by the proof of Friedman \cite{Friedman}, that the real and imaginary parts of the components of $\Pb $ are real analytic functions. Also, since we know that $\e^\dagger \Pb  \e$  is strictly greater than 0 on $\overline{\mathbb{C}}$, the real and imaginary parts of the function $(\e^\dagger \Pb  \e)^{-1}$ are also real analytic and so, by induction, the real and imaginary parts of the components of $D_{\pm}^j\Pb $ are analytic for any $j.$
Next, we need to integrate the real, analytic function $||D_{\pm}^j\Pb ||_{(\Pb )}^2$ over the complex plane.

Using Lemma \ref{hatP} we make use of the fact that, if $\Pb (\xi_{+},\xi_-)$ is a rank-1 projector solution which is in $C^{N+2}(\overline{\mathbb{C}})$ and with finite action then we can define a new function by inverting the variables across the unit circle. That is, if we define ${\hat \Pb } (w_+,w_-)$ as in \eqref{hatP}
then $\hat{\Pb }$ will also be a rank-1 projector solution in $C^{N+2}(\overline{\mathbb{C}})$ with finite action. 


Again, to simplify the notation, we will write $\Pb =\Pb (\xi_+, \xi_-), \ \hat{\Pb }=\hat{\Pb }(w+,w_-)$ and
$$ D_{\xi_{\pm}}^{(\Pb )}\equiv D_{\xi_{\pm}}, \qquad D_{w_{\pm}}^{(\hat{\Pb })}\equiv D_{w_{\pm}}$$
where it is understood that the covariant derivative in the direction of $\xi_\pm$ is taken with respect to $\Pb $ and with respect to $\hat{\Pb }$ for  derivatives in the direction $w_\pm.$ Next, we can show by induction that the inversion of $D_{\xi_{\pm}}^j \Pb $ about the disk is given by
\be \label{ind4} D_{\xi_{\pm}}^j \Pb =\sum_{k=1}^{j} c_kw_{\pm}^{k+j}D_{w_{\pm}}^{k}\hat{\Pb }, \quad  c_k\in \mathbb{C}.\ee
From simplicity, we define $c_k\equiv 0$ for $k=0$ and $k > j$
From the case $j=1$
$$D_{\xi_{\pm}} \Pb =-w_{\pm}^{2}D_{w_{\pm}}\hat{\Pb }$$
Now, assume it is true for $j$ and we prove for $j+1$
\bea D_{\xi_{\pm}}^{j+1} \Pb &=&D_{\xi_{\pm}}(D_{\xi_{\pm}}^{j} \Pb )\nn
&=&\sum_{k=1}^{j} -(k+j)c_kw_{\pm}^{k+j+1}D_{w_{\pm}}^{k}\hat{\Pb }-c_kw_{\pm}^{k+j+2}D_{w_{\pm}}^{k+1}\hat{\Pb } \nn
&=&\sum_{k=1}^{j+1}c_k'w_{\pm}^{j+1+k}D_{w_{\pm}}^{k}\hat{\Pb }\eea
where $c_k'=-(k+j)c_k-c_{k-1}.$ Thus, we have proven \eqref{ind4}.

Next, we use the fact that $|\xi_+|\geq1$ then $|w_+|\leq 1$ and \eqref{ind4} to transform the integral. 
$$\int_{|\xi_+|\geq 1}\!\!\!\!\!\! ||D_{\xi_\pm}^j\Pb ||_{(\Pb )}^2d\xi_+d\xi_-=\int_{|w_+|\leq 1}||\sum_{k=1}^{j} c_kw_{\pm}^{k+j}D_{w_{\pm}}^{k}\hat{\Pb }||_{(\hat{\Pb })}^2|w_+|^{-2}dw_+dw_-$$
Using the triangle inequality of the semi-norm, 
$$||\sum_{k=1}^{j} c_kw_{\pm}^{k+j}D_{w_{\pm}}^{k}\hat{\Pb }||_{(\hat{\Pb })}^2\leq \sum_{k=1}^{j} |c_k|^2 |w_{+}|^{2(k+j)}||D_{w_{\pm}}^{k}\hat{\Pb }||_{(\hat{\Pb })}^2$$
we obtain the inequality 
\be \label{extDnP}\int_{|\xi_+|\geq 1} \!\!\!\!\!\!||D_{\xi_\pm}^j\Pb ||_{(\Pb )}^2d\xi_+d\xi_-\leq \sum_{k=1}^{j}\int_{|w_+| \leq 1}  \!\!\!\!\!\!\!\!\!\!\!\!|c_k|^2 |w_{+}|^{2(k+j-1)}||D_{w_{\pm}}^{k}\hat{\Pb }||_{(\hat{\Pb })}^2dw_+dw_-.\ee
Since, since, $k+j>1$  the right-hand side of \eqref{extDnP} is an integral of an analytic function over a compact set  and so is finite. 
Hence, we see that we can break our integral up into a finite sum of integrals of real analytic functions over unit disks and so the integral must be finite.  
\bea \int_{\mathbb{R}^2} ||D_{\xi_{\pm}}^j\Pb ||_{(\Pb )}^2d\xi_+d\xi_-= \int_{|\xi_+|\leq 1} ||D_{\pm}^j\Pb ||_{(\Pb )}^2d\xi_+d\xi_- \nn
+\sum_{k=1}^j \int_{|w_+|\leq 1}|c_k|^2 |w_{+}|^{2(k+j-1)}||D_{w_{\pm}}^{k}\hat{\Pb }||_{(\hat{\Pb })}^2dw_+dw_-< \infty \eea

Thus  $|A_{0,n}|=|A_{i,j}|$ is in  $L^1(\mathbb{R}^2).$ But, we also know that 
\be \label{idk1} |\int_{\mathbb{R}^2}A_{0,n}d\xi_+d\xi_-| \leq \int_{\mathbb{R}^2}|A_{0,n}|d\xi_+d\xi_-< \infty \ee
and so, by Liouville's Theorem, the analytic function $A_{0,n}$ must be identically 0 and so by \eqref{indexshift} $A_{i,j}=0$ for all $i+j\geq 1.$
Finally, since 
$$A_{i,j}=\frac{(D_{-}^i\Pb  \e)^\dagger(D_{+}^j\Pb  \e)}{\e^\dagger \Pb  \e}\equiv 0 \quad \forall i+j\geq 1$$
 we see that the subspaces 
\be H=\lbrace \v| \v=\sum c_k D_{+}^k \Pb  \e, \quad 0\leq k \rbrace,\ee 
\be H^\perp =\lbrace \w| \w=\sum c_k D_{-}^k \Pb  \e, \quad 1 \leq k  \rbrace\ee 
are mutually orthogonal.  \QED 

\begin{lemma}\label{2}
Let $n+1$ be the dimension of $H$ and $m$ the dimension of $H^\perp$. We shall show that the projectors $\Pi_{-}^k\Pb , 0\leq k$ project in to $H$ and that the projectors $\Pi_{+}^k\Pb , \ 1\leq k$ project in to $H^\perp.$ Furthermore the images of $\Pi_{-}^k\Pb , 0\leq k\leq n$ form a basis for $H$ and $\Pi_{+}^k\Pb , 1\leq k\leq m$ form a basis for $H^\perp.$
\end{lemma}

Let $n+1$ be the dimension of $H$ and $m$ the dimension of $H^\perp$. Immediately, we can see that $m+n+1\leq N$ and when there is equality, $H\cup H^\perp = \C.$ Also, in particular $n, m <N.$
It is easy to see that the vectors $D_{-}^k\Pb  \e$ must be linearly independent for $k=0,...,n$ otherwise the dimension of $H$ would be strictly less than $n+1$ and similarly for 
$D_{+}^k\Pb  \e$  for $k=1,..., m$ and so these vectors form a basis for $H$ and $H^\perp.$

Next, we show that the projectors $\Pi_{-}^k\Pb , 0\leq k$ project in to $H$ and that the projectors $\Pi_{+}^k\Pb , \ 1\leq k$ project in to $H^\perp$ and that and a finite spanning set of the vectors onto which they project are linearly independent. Recall, we assumed that the projectors are not orthogonal to $\e$ and so we can also represent the vectors which the projectors map onto by multiplying them by $e.$

We prove that for any vector $\v$, there exist $ \psi_j=\psi_j(\xi_+, \xi_-)$ scalar functions such that
\be \label{Pik expand} (\Pi_{\pm}^k\Pb )\v =\sum_{j=1}^{k}\psi_jD_{\pm}^j\Pb  \e,\qquad  \psi_j:\mathbb{R}^2\rightarrow \mathbb{C},\ k\geq 1.\ee
and $\psi_k$ is identically $0$ for $\v=\e$ only when $\Pi_\pm^k\Pb=\0.$ 

Proof, From the case $k=0, \ \Pb $, we know that $\Pb $ maps any vector onto $\Pb  \e$ so $\Pb  \v$ is proportional to $ \Pb  \e \in H.$

From, $k=1,$ we have by  \eqref{rlcov}
\be (\Pi_{\pm}\Pb )\v=\frac{\e^\dagger (D_{\pm}\Pb )^\dagger \v}{\e^\dagger (D_{\pm}\Pb )^\dagger D_{\pm}\Pb  \e}D_{\pm}\Pb  \e.\ee
Assume true for $k$, using  \eqref{PX} we obtain
$$ (\Pi_{\pm}^{k+1}\Pb )\v= M(\v) (\partial_{\pm}\Pi_{\pm}^k\Pb )\Pi_{\pm}^k\Pb  \e$$
where $M(\v)$ is a scalar function depending on the vector $\v$ defined as
\be \label{Mv} M(\v)=\frac{\e^\dagger \Pi_{\pm}^k\Pb \partial _{\mp}\Pi_{\pm}^k\Pb  \v}{\e^\dagger \Pi_{\pm}^k\Pb  \partial _{\mp}(\Pi_{\pm}^k\Pb ) \partial_{\pm}(\Pi_{\pm}^k\Pb )\Pi_{\pm}^k\Pb  \e}.\ee
Notice that $M(\e)=0$ only when $\Pi_{\pm}^{k+1}\Pb =0.$
But also, we know that since $\Pi_{\pm}^k\Pb $ is a rank-1 projector, we can use  \eqref{DX} to show that, there exists a scalar function $\phi(\Pi_{\pm}^k\Pb )$ such that 
$$(\partial_{\pm}\Pi_{\pm}^k\Pb )\Pi_{\pm}^k\Pb  \e=(D_{\pm}\Pi_{\pm}^k\Pb -\phi(\Pi_{\pm}^k\Pb ) \Pi_{\pm}^k\Pb )\e.$$
So, finally, we have 
\bea (\Pi_{\pm}^{k+1}\Pb )\v&=&M(\v)\left(D_{\pm}\Pi_{\pm}^k\Pb -\phi(\Pi_{\pm}^k\Pb ) \Pi_{\pm}^k\Pb \e \right)\nn
&=&M(\v)\left(\sum_{j=1}^k \psi_jD_{\pm}^{j+1}\Pb  \e+\left(\partial_\pm \psi_j-\psi_j\phi(\Pi_{\pm}^k\Pb )\right)D_{\pm}^j\Pb  \e\right).\nonumber \eea
Here, we have use the induction assumption that \eqref{Pik expand} holds for $k.$
Thus, if we define new scalar function $\psi_k':\mathbb{R}^2\rightarrow \mathbb{C}$ for $j=1,..., k+1$ as 
\[\psi_j'=M(\v)\left(\psi_{j-1}+\left(\partial_\pm \psi_j-\psi_j\phi(\Pi_{\pm}^k\Pb )\right)\right), \qquad j\leq k\]
\[ \psi_{k+1}'=M(\v)\psi_{k}\]
 we see that $$(\Pi_{\pm}^{k+1}\Pb )\v=\sum_{j=1}^{k+1}\psi_j'D_{\pm}^j\Pb  \e$$ and that $\psi_{k+1}$ is identically $0$ for $\v=\e$ only when $\Pi_\pm^k \Pb =0$ and so we have proved \eqref{Pik expand} by induction. 
 
Thus, the projectors $\Pi_{-}^k\Pb $ for $k\geq 0$ maps in to $H$ and the projectors $\Pi_{-}^k\Pb $ for $k\geq 1$ map into $H^\perp.$ Furthermore,  the images of $\Pi_{-}^k\Pb , 0\leq k\leq n$ form a basis for $H$ and $\Pi_{+}^k\Pb , 1\leq k\leq m$ form a basis for $H^\dagger.$
 \QED 

\begin{lemma}\label{3} The projectors $\Pi_{\pm}^k\Pb , k\in \mathbb{N}$ are mutually orthogonal and we have a decomposition of the operator 
$$I_{n+m+1}=\Pi_{-}^{n}\Pb +...+\Pi_{-}\Pb +\Pb +\Pi_{+}\Pb +...+ \Pi_{+}^{m}\Pb $$ 
which acts as the identity on $H\cup H^\perp.$
\end{lemma}
Proof. Since the spaces $H$ and $H^\perp$ are mutually orthogonal by Lemma \ref{1}, we see that for any two vectors $v,w$ and $j+k>1$, we have $(\Pi_{-}^j\Pb  \v)^\dagger \Pi_{+}^k\Pb  \w=0$  and so 
 $ \Pi_{-}^j\Pb \Pi_{+}^k\Pb  =\0.$ Thus, it remains to show that 
$$\Pi_{-}^j\Pb \Pi_{-}^k\Pb  =\0, \qquad \Pi_{+}^j\Pb \Pi_{+}^k\Pb  =\0, \quad \forall j+k\ne0.$$
If we generalize the results, we have shown that for any finite action solution $\Pb '$ to the equations of motion, which is at least $M+2$ differentiable we can create a sequence of projectors so that $\Pb '\Pi_{\pm}^k\Pb '=\0$ for $1\leq k \leq M.$
Now, consider $\Pi_{-}^j\Pb \Pi_{-}^k\Pb  =\0$ with $j>k$ without loss of generality. 
Then, $\Pi_{-}^k\Pb =\Pi_{+}^{j-k}\Pi_{-}^j\Pb $ since by Theorem \ref{Pi e-l} each $\Pi_{\pm}^l\Pb $ satisfies the equations of motion. Also, from Theorem \ref{Pi action}, we see that $\Pi_{-}^j\Pb $ has finite action. By assumption, $\Pi_{-}^j\Pb $ is at least $N-j+2$ differentiable and the previous results hold for $\Pi_-^j\Pb  =\Pb '$ and we obtain 
$$\Pi_{-}^j\Pb \Pi_{-}^k\Pb =\Pi_{-}^j\Pb \Pi_{+}^{j-k}(\Pi_{-}^j\Pb )=\Pb '\Pi_{+}^{j-k}(\Pb ')=\0.$$
Similarly, for $j>k$
$$\Pi_{+}^j\Pb \Pi_{+}^k\Pb =\Pi_{+}^j\Pb \Pi_{-}^{j-k}(\Pi_{+}^j\Pb )=\0.$$
Thus, the the operators $ \Pi_{\pm}^k\Pb $ for  $0\leq k\leq N $ are mutually orthogonal. 
Finally, for any element $\v \in H\cup H^\dagger,$ $\v$ can be written as 
$$\v=\sum_{i=0}^{n}c_i \Pi_{-}^i \Pb  \e+\sum_{j=1}^{m}\Pi_{+}^j\Pb  \e$$ and so 
$${\bf I}_{m+n+1}\v=\left(\Pi_{-}^{n}\Pb +...+\Pi_{-}\Pb +\Pb +\Pi_{+}\Pb +...+ \Pi_{+}^{m}\Pb \right)\v=\v$$
and thus $ {\bf I}_{m+n+1}$ acts as the identity on $H\cup H^\dagger.$ In particular, if $\Pb $ is sufficiently arbitrary so that $H \cup H^\dagger$ spans $\C$ then ${\bf I}_{m+n+1}$ will be the identity matrix and we have a decomposition of unity by orthogonal projectors \QED

Proof of Theorem \ref{bigone}. Assume $\Pb $ is a finite action solution of the Euler-Lagrange equations $\Pb :\overline{\mathbb{C}} \rightarrow Gl_n(\mathbb{C})$ which is a rank-1 projector whose components are differentiable at least $N+2$ times.  From Lemma \ref{2} $\Pi^{n+1}_{-}\Pb $ maps any vector into $H$ so for arbitrary vectors $\v, \w$  $$\Pi^{n+1}_{-}\Pb  \w=I_{m+n} \Pi^{n+1}_{-}\Pb  \w$$
which implies
$$\v^\dagger \Pi^{n+1}_{-}\Pb  \w=\v^\dagger I_{m+n+1}\Pi^{n+1}_{-}\Pb  \w.$$
However,  from the results of Lemma \ref{3}, we know that $\Pi^{n+1}_{-}\Pb \Pi^{k}_{\pm}\Pb =\0$ for all $0 \leq k$ and so for any vectors $\v, \w$ $\v^\dagger \Pi^{n+1}_{-}\Pb  \w=0$ and so $ \Pi^{n+1}_{-}\Pb  \equiv \0.$ Thus, $\Pi^{n}_-\Pb $ is holomorphic. 
Similarly, we see that $\Pi^{m+1}_{+}\Pb =\0$ and so $\Pi^{m}_{+}\Pb $ is antiholomorphic.

Define 
\be \Pb _0=\Pi^{n}_{-}\Pb,  \qquad \Pb _j=\Pi_+^j \Pb_0 .\ee
Since each monomial $\Pi_{\pm}^k\Pb $ satisfies the equations of motion, we can use Euler-Lagrange \eqref{el3} to show that $\Pb =\Pb _n $ and can be written as either 
\be \Pb =\Pi^{n}_{+}\Pb _0, \qquad \Pb = \Pi^m_-\Pb _{m+n}.\ee
\QED

Thus,  for any solution to the Euler-Lagrange equation $\Pb $ in $C^{N+2}(\overline{\mathbb{C}})$ for $ 0\leq M\in \mathbb{N}$ with finite action, there exists a holomorphic projector which generates a sequence of rank-1 Hermitian projectors for $m+n+1\leq N$ 
\be\label{lambda} \Lambda \equiv \lbrace \Pb _0, \Pi_+\Pb _0, \Pi_{+}^2\Pb _0,..., \Pi_+^{m+n-1}\Pb _0, \Pi_+^{m+n}\Pb _0\rbrace\ee
which are all finite action solutions to the Euler-Lagrange equations, are mutually orthogonal and for which $\Pi_{+}^n\Pb _0=\Pb $ and  $\Pi_-P_0=\Pi_+^{m+n+1}P_0=\0.$
Furthermore, if  $\partial_\pm^k \Pb $ spans $\C,$ then the set $\Lambda$ will generate a basis for $\C.$ 

\section{Surfaces associated with $\CP$ sigma models}
Beginning with the sine-Gordon model, there has been much interest in surfaces associated with $\CP$ models and their geometry. The generalized Weierstrauss formula for immersion  developed by by B Konopenchenko and I Taimanov \cite{Kono1996, KonoTaim1996} provides an efficient tool for constructing a surface immersed in $su(N)$ for any solution of the $\CP$ sigma model. We consider such surfaces in this section. 

We can see from the Euler-Lagrange \eqref{el2} written as a conservation law that the following differential
\be \label{dF} d\F=i \left(-[\partial_+\Pb ,\Pb ]d \xi_++[\partial_-\Pb ,\Pb ]d\xi_-\right)\ee 
is closed and skew-Hermitian. The differential can then be integrated to obtain a surface immersed in the  $su(N)$ algebra associated to $\Pb ,$ a solution of the $\CP$ sigma model. We have following theorem, first proven in \cite{GSZ2005}.
\begin{theorem} \label{Fkthm} If $\Pb $ is a rank-1 Hermitian projector is a solution to the Euler-Lagrange equations then there exists a surface $\F$ whose immesion is defined by 
\be \label{F} \F=i\int -[\partial_+\Pb ,\Pb ]d \xi_++[\partial_-\Pb ,\Pb ]d\xi_- , \qquad \F^\dagger =-\F\in su(N).\ee 
\end{theorem} 
Notice immediately that because the immesion of the surface is give in terms of the projectors, it will necessarily be gauge invariant. 

For uniformity, we introduce the following inner product on $su(N),$ proportional to the Killing form, and extend it to a bilinear product on all of $Gl_N(\mathbb{C}(\xi_+, \xi_-))$ as
\be \label{()} ({\bf A},{\bf B})=-\frac 12 tr({\bf A}\cdot {\bf B}), \qquad {\bf A},{\bf B} \in su(N) \ee
The following theorem was first proven in \cite{GrunYurd2009}.
\begin{theorem} \label{surface} From $\Pb $ a finite action solution on $\overline{\mathbb{C}}$ of the Euler-Lagrange equations \eqref{EL}, the surface $\F$ defined by \eqref{F} is conformally parameterized. 
\end{theorem}
Proof. We can compute the components of the metric tensor using
\be \partial_+\F=-i[\partial_+\Pb ,\Pb ] \qquad \bp \F=i[\partial_-\Pb ,\Pb ]\ee
to obtain, 
$$g_{++}=(\partial_+ \F, \partial_+ \F)=\frac12 tr([\partial_+\Pb , \Pb ]\cdot [\partial_+\Pb , \Pb ])=-\langle D_-\Pb ,D_+\Pb \rangle =0$$
$$g_{+-}= (\partial_+ \F, \partial_- \F)=\frac12tr([\partial_+\Pb , \Pb ]\cdot [\partial_+\Pb , \Pb ])=\frac12 tr(\partial_+\Pb \partial_-\Pb ),$$
and the complex conjugate of $g_{++},$ $\overline{g_{++}}=g_{--}=0.$ 
\QED

\subsection{Geometry of the surfaces in terms of the projectors}
In this section, we discuss the geometry of the surfaces using the projector formalism created in the Section 2 as well as results obtained from the proofs in Section 4. Unlike the previous work on the surfaces \cite{GoldGrund2010}, we give the real form of the geometry, in terms of the real coordinates $\xi_1$ and $\xi_2 $ with $\partial_i=\partial /\partial \xi_i.$

The tangent vectors are then 
\be \partial_1 \F=[\partial_2\Pb ,\Pb ], \qquad \partial_2 \F=-[\partial_1 \Pb ,\Pb ]\ee 
note that $\partial_i \F \in su(N).$ 
We use the fact that we can rewrite 
\be L(\Pb ) =||D_+\Pb ||^2+||D_-\Pb ||^2=||D_+\Pb +D_-\Pb ||^2\ee
to expand $L(\Pb )$ as a real form using the fact the $\langle D_+\Pb , D_-\Pb \rangle=0.$
\be L(\Pb )=||D_1\Pb ||^2=||D_2\Pb ||^2\ee
here we use definition of the covariant derivative \eqref{D} to write, 
$$D_i=\partial_i+\langle \Pb , \partial_i\Pb \rangle , \qquad D_1=D_++D_-, \quad D_2=i(D_+-D_-).$$
We can then use the inner product given above \eqref{()} to determine the real form of the metric tensor as 
\be g_{11}=(\partial_1 \F, \partial_1\F)= ||D_1\Pb ||^2,\quad  g_{22}=(\partial_2 \F, \partial_2\F)= ||D_1\Pb ||^2,\ee
\be g_{12}=(\partial_1 \F, \partial_2\F)=0 =g_{21} \ee
and so it is conformally parameterized with first fundamental form 
\be I =||D_1\Pb ||^2\left(d\xi_1^2+d\xi_2^2\right).\ee
The  Christoffel symbols are given by 
\be \Gamma^1_{11}=\Gamma^2_{12}=-\Gamma^{1}_{22}=\frac{\partial_1 ||D_1\Pb ||^2}{2 ||D_1\Pb ||^2}, \quad  \Gamma^2_{22}=\Gamma^1_{12}=-\Gamma^{2}_{11}=\frac{\partial_2  ||D_1\Pb ||^2}{2 ||D_1\Pb ||^2}\ee
where 
\bea \partial_1  ||D_1\Pb ||^2&=&\langle D_+\Pb , D^2_+\Pb \rangle+ \langle D_-^2\Pb , D_- \Pb \rangle  
 +\langle D_-\Pb , D_-^2\Pb \rangle +\langle D_+^2\Pb , D_+\Pb \rangle \nn
 \partial_2  ||D_1\Pb ||^2&=&i(\langle D_+\Pb , D^2_+\Pb \rangle+ \langle D_-^2\Pb , D_-\Pb \rangle  
 -\langle D_-\Pb , D_-^2\Pb \rangle -\langle D_+^2\Pb , D_+\Pb \rangle) \nn
 \frac{\partial_1  ||D_1\Pb ||^2}{2 ||D_1\Pb ||^2}&=& \frac{Re(\langle D_1\Pb , D_2^2\Pb \rangle) }{||D_1\Pb ||^2}, \qquad  \frac{\partial_2 ||D_1\Pb ||^2}{2||D_1\Pb ||^2}= \frac{Im(\langle D_2\Pb , D_1^2\Pb \rangle) }{||D_1\Pb ||^2}.\eea
 Here we used repeatedly the fact that the $A_{i,j}=0$ from Lemma \ref{1}. 

 We use the Riemannian tensor $R(\X,\Y){\bf Z}$ \cite{DoCarmo} to compute the Gaussian curvature of the surface using the basis for the tangent space given by 
 $\partial_1\F$ and $\partial_2 \F.$
 \bea K(\F)&=& \frac{(R(\partial_1\F, \partial_1\F)\partial_1\F, \partial_1\F)}{g_{11}}\eea
$$=4\left(1-3\frac{||D_1\Pb +D_2\Pb ||^2||D_1\Pb -D_2\Pb ||^2}{||D_1||^4}\right)-2\frac{||D_1D_2\Pb ||^2||D_1\Pb ||^2-|\langle D_1, D_1^2\Pb \rangle|^2}{||D_1\Pb ||^6}.$$
 We can consider explicitly the embedding of our surface into $su(N)$ and define a symmetric bilinear mapping $\mathcal{B} $ from the tangent space of $\F$ to its normal given in terms of the basis functions as
 \be \mathcal{B} (\partial_i\F,\partial_j \F)=\partial_i\partial_j\F-\sum_{k=1,2} \Gamma_{ij}^k\partial_k\F\equiv \N _{ij},\ee
 where $\N _{i,j}$ are by construction orthogonal to the surface. 
This symmetric bilinear mapping $B$ can then be used to construct the second fundamental form 
 \be II=\frac{1}{||D_1\Pb ||^2}\left(\N _{11}d\xi_1d\xi_1+\N _{12}d\xi_1d\xi_2+\N _{22}d\xi_2d\xi_2\right).\ee
 The mean curvature vector of the immersion is given by the trace of the second fundamental form, 
\be \H=\frac{1}{||D_1\Pb ||^2}\left(\partial_1^2\F+\partial_2^2\F\right) \in su(N). \ee
We use these quantities to compute the Willmore functional and the Euler-Poincare characteristic of the surface
\be W(\F)= \frac{1}{2}\int_{\mathbb{R}^2} (H,H) g_{ij}d\xi^i d\xi^j \qquad  \Delta(\F) =\frac{1}{2\pi} \int_{\mathbb{R}^2} K(\F)g_{ij}d\xi^i d\xi^j.\ee

\subsection{Geometric characteristics in terms of the physical quanties.}
A surprising result is that many of the geometric characteristics of the surface, such as area, Christoffel symbols, Willmore functional and Euler-Poincare character can be obtained directly from the physics quantities such as the Lagrangian density, topological charge density and the action. We present these results in this section and for simplicity of calculations return to complex coordinates. 

\begin{theorem}
 For a finite action solution to the Euler-Lagrange equations \eqref{EL}, $\Pb $, the conformal factor associated with the surface $\F$ defined by \eqref{F} is proportional to the Lagrangian density and area of the surface is given by the action functional of the model. In particular, the surface will have finite area.  
\end{theorem}
Proof. From Theorem \ref{surface}, we can see immediately that the surface is conformally parameterized and the factor is given by  the Lagrangian density  $L(\Pb)$ \eqref{lagrangian1} 
\[g_{+-}=\frac12 tr(\partial_+\Pb \partial_-\Pb )=\frac 12 L(\Pb )\]
Thus, the first fundamental form is given by 
\be \label{FFF} I(\F)=L(\Pb )d\xi_+d\xi_-\ee
and the area of the surface is equal to the action $S(\Pb)$\eqref{action}, 
\be \label{area} A(\F)=\int_ {\mathbb{R}^2}2g_{+-}d\xi_+d\xi_-=\int_ {\mathbb{R}^2}L(\Pb )d\xi_+d\xi_-=S(\Pb ).\ee
Note that the area of the surface is a definite integral over the plane $\mathbb{R}^2$ and so does not depend on $\xi_\pm.$ 
In complex coordinates, the only non-zero Christoffel symbols can be expressed terms of the Lagrangian density as  
\be\label{christ} (\Gamma(\F))^+_{++}=\frac{\partial_+ L(\Pb )}{L(\Pb )}, \qquad  (\Gamma(\F))^-_{--}=\frac{\partial_- L(\Pb )}{L(\Pb )} .\ee
Thus, the covariant derivatives are given by 
\be \nabla_{\partial_{\pm} \F}\partial_{\mp}\F=0, \qquad \nabla_{\partial_\pm \F}\partial_{\pm}\F=\ (\Gamma(\F))^{\pm}_{\pm \pm }\partial_{\pm}\F\ee
We see immediately that the surface is torsion free; that is $$T^a_{bc}=\Gamma^a_{bc}-\Gamma^a_{cb}=0.$$
In complex coordinates, the Gaussian curvature becomes 
\be \label{KF} K(\F)=-\frac{2\partial_+\partial_-(ln(L(\Pb )))}{L(\Pb )}.\ee
The mean curvature vector, written in matrix form, is given by 
\be \label{mean} \H(\F)=\frac{-4i}{L(\Pb )}[\partial_+\Pb ,\partial_-\Pb ]\ee
and thus is traceless. 
The norm of this vector can be written in terms of  the Lagrangian density and the topological charge density $q(\Pb)$ \eqref{chargedensity},
\bea (\H(\F),\H(\F))&=&\frac{8}{L(\Pb )^2}tr([\partial_+\Pb ,\partial_-\Pb ]^2)=\frac{4}{L(\Pb )^2}\left(L(\Pb )^2+ 3q(\Pb )^2\right).\eea
From this quantity, we can rewrite the Willmore functional, in terms of the action, the Lagrangian density, and the topological charge density  \eqref{chargedensity},
\be \label{W} W(\F)=\frac{1}2S(\Pb )+\frac{3}{2}\int_{\mathbb{R}^2}\frac{q(\Pb )^2}{L(\Pb )}d\xi_+d\xi_-,\ee
and the Euler-Poincare Character in terms of just the Lagrangian density as
\bea \label{delta} \Delta(\F) &=&\frac{1}{2\pi} \int_{\mathbb{R}^2} K(\F)g_{12}d\xi_+d\xi_- = -\frac{1}{\pi} \int_{\mathbb{R}^2}\p \bp ln(L(\Pb ))d\xi_1d\xi_2.\eea
 
 Thus, we see that it is possible to represent several geometric quanties of the surface from the physical quanties of the model. On the other hand, one can immediately compute the action and the Lagrangian density of the physical system via the metric tensor and the surface area of the surface, respectively. We also compute the absolute value of the topological charge density via the relation
\be |q(\Pb )|=g_{12}\sqrt{\frac{(\H(\F),\H(\F))-4}{3}}. \ee

\subsection{Normals to the surface}
In this section, we consider normals to the surface again using complex coordinates to simplify the computations. The second fundamental form is
\be\label{II} II=\N _{++}d\xi_+d\xi_++\N _{+-}d\xi_+d\xi_-+\N _{--}d\xi_-d\xi_-\ee
where the 3 normals are give by
$$\N _{++}=\p \p \F-(\Gamma(\F))^+_{++}\p \F$$
$$\N _{+-}=\p \bp \F=\H$$
$$\N _{--}=\bp \bp \F-(\Gamma(\F))^-_{--}\p \F$$
and have real forms in the $su(N)$ algebra. 
Here we show that these normals are mutually orthogonal and in fact orthogonal to the surface. 
To do this,  we use results from Lemma \ref{1} which give, for $\Pb$ a finite action solution of the $\CP$ model,
\[ \langle D_+\Pb , D_-\Pb \rangle =0 \Rightarrow tr(\Pb \p \Pb  \p \Pb )=0\]
\[\langle D_-\Pb , D_+\Pb \rangle=0 \Rightarrow tr(\Pb \bp \Pb \bp \Pb ) =0\] 
\[ \langle D_-\Pb , D_+^2\Pb \rangle =0 \Rightarrow tr(\Pb \p \Pb  \p^2\Pb )=0\]
\[ \langle D_-^2\Pb , D_+\Pb \rangle=0 \Rightarrow tr(\Pb  \p^2\Pb \p \Pb )=0.\]

We verify that these normals are in fact orthogonal to the tangent vectors of the  surface 
\[ (\N _{++}, \p \F)=-\frac12 tr( [\p^2\Pb ,\Pb ]\cdot[\p \Pb ,\Pb ]) =0\]
 \[(\N _{--}, \bp \F)=-\frac 12 tr( [\p^2\Pb ,\Pb ]\cdot[\p \Pb ,\Pb ]) =0\]
\[(\N _{++},\bp \F)=-\frac12\left(tr( [\p^2\Pb ,\Pb ]\cdot[\bp \Pb ,\Pb ]) -\partial_+ \left(tr(\Pb \bp \Pb  \p \Pb )+tr(\Pb  \p \Pb \bp \Pb )\right)\right)\]
\[=0\hspace*{220pt}\]
\[(\N _{--},\p \F)=-\frac12\left(tr([\p^2\Pb ,\Pb ]\cdot[\bp \Pb ,\Pb ]) -\partial_+tr(\Pb \bp \Pb  \p \Pb )+\partial_+tr(\Pb  \p \Pb \bp \Pb )\right)\]
\[=0\hspace*{220pt}\]
 From the normal $\N _{+-}=\H$ the mean curvature \eqref{mean}, we can simplify the computations by expanding  in terms of the raising and lowering operators as
\[\N _{+-}=\frac{-4i}{L(\Pb )} \left(\Pi_+\Pb  ||D_+\Pb ||^2-\Pi_- \Pb  ||D_-\Pb ||^2 -q(\Pb ) \Pb \right) \]
and  it is orthogonal to the tangent vectors since 
\[ tr(\Pi_\pm \Pb  [\partial_\pm \Pb ,\Pb ])=0\qquad tr(\Pb  [\partial_\pm \Pb ,\Pb ])=0.\]
In fact, these equalities are true for any projector in the sequence constructed in Theorem 8. That is, any element of the set $\Lambda$ given by \eqref{lambda} will be normal to the surface since 
\be tr(\Pi_\pm ^k\Pb  [\partial_\pm \Pb ,\Pb ])=tr(\Pb \Pi_\pm ^k\Pb  \partial_\pm \Pb )-tr(\Pi_\pm ^k\Pb \Pb  \partial_\pm \Pb )=0.\ee

In addition,  the normals $\N _{++}$ and $\N _{--}$ are orthogonal to each element in $\Lambda.$

\be tr([\partial_\pm^2\Pb ,\Pb ] \Pi_\pm^k\Pb )=tr( \partial_\pm^2\Pb  \Pi_\pm^k\Pb ) -tr( \partial_\pm^2\Pi_\pm^k\Pb \Pb )=0.\ee

Thus, for maximally non-degenerate projector $\Pb $ where $n+m+1=N$ we have $N$ normals given by the sequence of projectors $\Pb _k$ and also two other normals given by the second fundamental form. Of course, the last two may be identically $0$ as would necessarily be the case for $\mathbb{C}P ^1$ 
where the surface is immersed in $\mathbb{R}^3$ and so there is a unique non-zero normal given by $\Pb .$ 

In fact, we not only have a sequence of normals but also a sequence of $n+m+1$ surfaces all of which have the same $n+m+1$ normals. From each  $\Pb _k$ we can associate it with a surface $\F_k$ which will be orthogonal to each of the  $\Pb _i$'s. In the following section we exploit the sequence of projectors to obtain an explicit expansions of the surfaces in terms of projectors.

\begin{subsection}{Characterization of all surfaces associated to $\CP$ sigma models}
We can use the results from Section 4 to completely determine all surfaces associated with $\CP$ sigma models defined on the extended complex plane  with finite action. Immediate from these results, we show that  any surface associated with a $\CP$ model will belong to a family of at most $N-2$ linearly independent surfaces built recursively from a ``holomorphic'' or ``antiholomorphic'' surface. By ``holomorphic'' surface, we mean that the surface is constructed from a set of $N-1$ holomorphic functions and their conjugates. However, the resulting surfaces still depends on both $\xi_+$ and $\xi_-.$ In this section we assume that $\Pb $ is sufficiently arbitary so that $m+n+1=N$ and $H\cup H^\dagger=\C$ and hence the elements of $\Lambda$  decompose $\C$ via projectors.  Otherwise, we could embed the entire model into a lower dimensional complex projective space and the following results would still hold for the lower dimension.

\begin{theorem} (Gundland and Yurdusen \cite{GrunYurd2009})
The surface defined as in Theorem 7 has an associated integer $k$ and holomorphic projector $\Pb _0$ so that the following holds
\be \label{Fk} \F=\F_k\equiv -i(\Pb _k+2\sum_{j=0}^{k-1}\Pb _j-\frac{1+2k}{N}{\bf I}_{N})), \qquad \Pb _\ell=\Pi^\ell_+\Pb _0.\ee 
\end{theorem}
We know from the theorems of the previous section, that for each $\Pb $ a rank-1 Hermitian projector which is a solution to the Euler-Lagrange equations with finite action on the extended complex plane there exists an integer $k$ and holomorphic projector $\Pb _0$ so that $\Pb =\Pi^k_+\Pb _0.$ It has been shown in the literature \cite{GrunYurd2009, GoldGrund2010} that the quantity 
$\F_k$ \eqref{Fk} has the required differential. 

Since the set of vectors $\Pb _i\e$ form an orthogonal basis for $\C$ we can immediately determine the characteristic polynomial for $\F_k$ as 
\be p_k(t)=(t-ic_k)^{N-k-1}(t-i(c_k-1))(t-i(c_k-2))^{k}, \qquad c_k=\frac{1+2k}{N},\ee 
we can use this to determine the minimal polynomial for $\F_k$ which is degree 2 for holomorphic ($k=0$) and antiholomorphic ($k=N-1$) and degree 3 for mixed solutions. 
The minimal polynomial for $k=0$ is 
\be \label{q0} q_0(t)=(t-ic_0)^2-i(t-ic_0), \qquad q_0(\F_0)=0,\ee
and for $k=1,.., N-2$
\be \label{qk} q_k(t)=(t-i(c_k-1))^3+(t-i(c_k-1)), \qquad q_k(\F_k)=0,\ee
for $k=N-1$
\be \label{qN-1} q_{N-1}(t)=(t-i(c_{N-1}-2))^2+i(t-i(c_{N-1}-2)), \qquad q_{N-1}(\F_{N-1})=0\ee
We can  see by direct calculation that $\F_0$ (and by analogy $\F_{N-1}$) will also satisfy the third order equation 
$(t-i(c_0-1))^3+(t-i(c_0-1))$ though this is not minimal. \\
We can invert the equation for $\F_k$ in terms of $\Pb _k$ as 
\be \label{PkF} \Pb _k=(\F_k-i(c_k-1))^2+I_N\ee
and use this to rewrite the differential of $\F_k$ and obtain the following theorem. 
\begin{theorem} $\F$ is a differentiable surface associated to a finite action solution of the $\CP$ sigma models on $\overline{\mathbb{C}}$ if and only if there exists some $k=0,.., N-1$ and $\lambda \in \mathbb{R}$ such that 
\be \label{dtF} det(\F)=(i)^N(\lambda)^{N-k-1}(\lambda-1)(\lambda-2)^{k}\ee
and the differential of $\F$ is given by 
\be \label{dFlambda} d\F=-i\left[\partial_+(\F-i(\lambda-1))^2,(\F-i(\lambda-1))^2\right]d\xi_+ \ee
$$+i\left[\partial_-(\F-i(\lambda-1))^2,(\F-i(\lambda-1))^2\right]d\xi_-.$$
\end{theorem}
Proof. We immediately have one direction of the proof from Theorem \ref{Fkthm}, with $\lambda=c_k.$ From the other direction,  assume \eqref{dtF} and \eqref{dFlambda}.
Now, \eqref{dtF} implies that $\F$ satisfies the polynomial 
\be \label{Flmin} (\F-i(\lambda-1))^3+(\F-i(\lambda-1))=0\ee
and so we can define $\Pb =(\F-i(\lambda-1))^2+I_N.$  The polynomial identity \eqref{Flmin} implies that $\Pb ^2=\Pb .$ From the skew-Hermitian property of $\F,$ we see that the projector is Hermitian. From the fact that the operator $\F$ has a unique eigenvector with eigenvalue $\lambda-1$ we see that $\Pb $ is of rank-1. Finally, from the fact that the differential of $\F$ must be closed, we obtain the Euler-Lagrange equations for $\Pb $ and so we see that the surface $\F$ is associated with a solution of the $\CP$ model. From the final statement, we see immediately that if the surface is conformally parameterized the conformal factor will be the Lagrangian density and so the surfaces area will be equal to the action of the associated sigma model.   \qed

We finish this section by noting an important algebraic identity that follows directly from the decomposition of the surfaces $\F_k$ in terms of the projectors $\Pb _j.$ We have the following identity. 
\be \sum_{k=0}^{N-1} (-1)^k\F_{k}=0\ee
This identity can be verified by using the fact that $ \sum_{j=0}^{N-1}c_j $ is $-1$ for even $N$ and $1$ for odd $N.$ Thus, there exists at most $N-1$ linearly independent surfaces $F_j, \ j=0,..., N-1$ immersed in $su(N).$
\end{subsection}
\begin{subsection}{(Anti-)Holomorphic Solutions}
 We now focus on the holomorphic and anti-holomorphic solutions of the model. That is, assume $\Pb _0$ is holomorphic and thus 
\be \label{F0}\F_0\equiv-i\left(\Pb _0-c_0{\bf I}_N\right), \qquad c_0=\frac 1N\ee
From such a surface, we have the two determining equations: the minimal polynomial
\be \label{minP0} q_0(t)=(t-ic_0)^2-i(t-ic_0), \qquad q_0(\F_0)=0\ee
and the holomorphicity requirement for $\Pb _0$
\be \label{F0holo} (\partial_-\Pb _0)\Pb _0=0 \iff [\F_0,\partial_-\F_0]+i\partial_-\F_0 =\0.\ee
From anti-holomorphic solutions, we use the fact that $c_{N-1}=2-c_0$ to obtain the two determining equations: the minimal polynomial
\be \label{minPN-1} q_{N-1}(t)=(t+ic_0)^2+i(t+ic_0), \qquad q_{N-1}(\F_{N-1})=0\ee
and the anti-holomorphicity requirement for $\Pb _{N-1}$
\be \label{FN-1antiholo} (\partial_+\Pb _{N-1})\Pb _{N-1}=0 \iff [\F_{N-1},\partial_+\F_{N-1}]+i\partial_+\F_{N-1}=\0.\ee
From (anti-)holomorphic solutions, we can greatly simplify the geometric quanties described above by noting that $$D_-^k\Pb _0\e=0, \qquad D_+^k\Pb _{N-1}\e=0$$ and so we have the Lagrangian density and nonzero component of the metric as
\be L(\Pb _0)=||D_+\Pb _0||^2=q(\Pb _0),\ee
which implies  \be g_{+-}(\F_0)=\frac 12 ||D_+\Pb _0||^2. \ee
For the anti-holomorphic surface we have, 
\be L(\Pb _{N-1})=||D_-\Pb _{N-1}||^2=-q(\Pb _{N-1}),\ee
which gives 
\be g_{+-}(\F_{N-1})=\frac 12 ||D_-\Pb _{N-1}||^2.\ee
We compute the other geometric quantities such as the Christoffel symbols, \eqref{christ}
\be (\Gamma(\F_0))^+_{++}=\frac{\langle D_+\Pb _0, D_+^2\Pb _0\rangle}{||D_+\Pb _0||^2 },\qquad (\Gamma(\F_0))^-_{--}= \overline{(\Gamma(\F_0))^+_{++}}.\ee
\be (\Gamma(\F_{N-1}))^+_{++}=\frac{\langle D_-\Pb _{N-1}, D_-^2\Pb _{N-1}\rangle}{||D_-\Pb _{N-1}||^2 },\qquad (\Gamma(\F_{N-1}))^-_{--}= \overline{(\Gamma(\F_{N-1}))^1_{11}}.\ee
The Gaussian curvature \eqref{KF} is given by 
\be K(\F_0)=4-2\frac{||D_+^2\Pb _0||^2||D_+\Pb _0||^2-|\langle D_+\Pb _0, D_+^2\Pb _0\rangle|^2}{||D_+\Pb _0||^6} \ee
\be K(\F_{N-1})=4-2\frac{||D_-^2\Pb _{N-1}||^2||D_-\Pb _{N-1}||^2-|\langle D_-\Pb _{N-1}, D_-^2\Pb _{N-1}\rangle|^2}{||D_-\Pb _{N-1}||^6} \ee
The Gaussian curvature is not necessarily either constant nor positive as seen in the examples of the next section. 

However, for any holomorphic or anti-holomorphic surface, the norm of the mean curvature \eqref{mean} is constant 
\be (\H(\F_0),\H(\F_0))= (\H(\F_{N-1}),\H(\F_{N-1}))=16.\ee
and so the  Willmore functional \eqref{W} is proportional to the action of the model
\be\label{Willholo} W(\F_0)=2S(\Pb _0), \qquad W(\F_{N-1})=2S(\Pb _{N-1}).\ee
Finally, we note that for (anti-)holomorphic solutions the topological charge \eqref{charge} is proportional to the action
\be\label{Chargeholo} Q(\F_0)=\frac 1 {\pi}S(\Pb _0)\qquad Q(\F_{N-1})=-\frac 1 {\pi}S(\Pb _{N-1}).\ee
\end{subsection} 
\begin{subsection}{Examples}
The above theoretical considerations are illustrated by the following examples of surfaces associated with holomorphic projectors. The first surface is give in terms of the Veronese sequence, the simplest solution to the $\CP$ sigma model. It is well known that the surface has constant positive Gaussian curvature. On the other hand, we introduce a new surface also defined by a holomorphic projector but with Gaussian curvature which is non-constant with varying sign. Both surfaces have finite surface area of $2\pi.$

\begin{subsubsection}{Veronese Surface in $\mathbb{C}\Pb ^2$}
We make a specific choice of holomorphic projector given by  
\be \Pb _0=\frac{1}{(1+|\xi|^2)^2} \left[\begin{array} {ccc} 
1 & \sqrt{2} \xi_- & \xi_-^2 \\
\sqrt{2}\xi_+ & 2\xi_+\xi_- & \sqrt{2}\xi_+\xi_-^2 \\
\xi_+^2 & \sqrt{2} \xi_+^2\xi_- & \xi_+^2\xi_-^2 \\                                
                               \end{array}\right]. \ee
 As in \eqref{F0}, the surface is given by $\F_0\equiv-i\left(\Pb _0-\frac{1}3{\bf I}_3\right)$.
The covariant derivatives of $\Pb _0$ are, for $e=[1,0,0]^\dagger,$
\be D_+\Pb _0=\frac{1}{(1+|\xi|^2)^3}\left[\begin{array}{ccc} 
-2\xi_-& -2\sqrt{2}2\xi_-^2&-2\xi_-^3\\
\sqrt{2}(1-|\xi|^2) & 2\xi_-(1-|\xi|^2) &\sqrt{2}\xi_-^2(1-|\xi|^2)\\
2\xi_+ & 2\sqrt{2}|\xi|^2 &2\xi_-^2\xi_+ \end{array} \right] \ee

 \be D_-\Pb _0=\frac{1}{(1+|\xi|^2)^2}\left[\begin{array}{ccc} 
0& \sqrt{2}&-2\xi_-\\
0& 2\xi_+ &\sqrt{2}2|\xi|^2\\
0& \sqrt{2}\xi_+^2&2\xi_+^2\xi_- \end{array}\right], 
\quad D_{-}\Pb _0e=0\ee

\[ |\langle D_+\Pb _0, D_+^2\Pb _0\rangle |^2=\frac{16 |\xi|^2}{(1+|\xi|^2)^6}\]
\[||D_+\Pb _0||^2=\frac{2}{(1+|\xi|^2)^2}, \qquad ||D_+^2\Pb _0||^2=\frac{4(2|\xi|^2+1)}{(1+|\xi|^2)^4}
\]

The metric, Gaussian curvature, action (and area of the surface \eqref{area}), Wilmore functional and topoligical charge are give respectively as,
$$g_{+-}(\F_0)=\frac{1}{(1+|\xi|^2)^2}$$
$$ K(\F_0)=4-2\frac{||D_+^2\Pb _0||^2||D_+\Pb _0||^2-|\langle D_+\Pb _0, D_+^2\Pb _0\rangle|^2}{||D_+\Pb _0||^6}=2$$
$$A(\F_0)=\int_{\mathbb{R}^2}\frac{2}{(1+|\xi|^2)^2}d\xi_+d\xi_-=2\pi$$

$$ W(\F_0)=4\pi, \qquad  Q(\F_0)=2, \qquad \Delta(\F_0)=2.$$
Thus, the surface associated with the holomorphic projector given by the Veronese sequence is compact, with finite surface area, constant positive curvature and hence constant positive Euler-Poincare characteristic and so is homeomorphic to the sphere. 
\end{subsubsection}
\subsubsection{Negative Curvature Example}
In this example we give an example of a surface associated with a holomorphic projector with non-constant curvature which has areas of both positive and negative value and whose asymptotic limit is negative. We also compute the area of the surface and the Euler-Poincare characteristic of the surface.  
The holomorphic projector we chose is given by 
\be \Pb _0=\frac{1}{1+|\xi_+|^2+4|\xi_+|^4}\left[ \ba{ccc} 
1 & \xi_- & 2\xi_-^2\\ 
\xi_+ & \xi_+\xi_- & 2\xi_-^2\xi_+\\
2\xi_+^2 & 2\xi_+^2\xi_-&4(\xi_+\xi_-)^2 \ea \right].\ee
As in \eqref{F0}, the surface is given by $\F_0\equiv-i\left(\Pb _0-\frac{1}3{\bf I}_3\right)$.
The metric tensor is given by
\[g_{+-}(\F_0)=\frac12 L(\Pb _0)=\frac{4|\xi_+|^4+16|\xi_+|^2+1}{2(1+|\xi_+|^2+4|\xi_+|^4)^2}
\]
we can use this to find the area of the surface
\[ A(\F_0)=\int_ {\mathbb{R}^2}\frac{4|\xi_+|^4+16|\xi_+|^2+1}{(1+|\xi_+|^2+4|\xi_+|^4)^2}d\xi_+d\xi_-=2\pi,\]
which is equal to the action of the model \eqref{action}.
The Gaussian curvature of the surface is given by
\[ K(\F_0)=-4\frac{448|\xi_+|^{12}-384|\xi_+|^{10}-2640|\xi_+|^8-4280|\xi_+|^6-660|\xi_+|^4-24|\xi_+|^2+7}{(4|\xi_+|^4+16|\xi_+|^2+1)^3}\]
which is negative at the origin and has asymptotic limit $K\rightarrow -14$ as $|\xi_+|\rightarrow \infty.$ 

It is interesting to note that the Euler-Poincar\'e character given by 
\[ \Delta(\F_0) =\frac{1}{2\pi} \int_{\mathbb{R}^2} K(\F)g_{12}d\xi_+d\xi_-=2\]
 is the same as surface generated by the Veronese sequence. 

The Willmore functional and topological charge are given immediately from \eqref{Willholo} and \eqref{Chargeholo}
$$ W(\F_0)=4\pi, \qquad  Q(\F_0)=2.$$
Hence, we see that even in the case of surfaces associated with holomorphic projectors, we can obtain surfaces with non-constant Gaussian curvature. 
\end{subsection}

\begin{section}{Concluding remarks}

We have shown that for any solution of $\CP$ sigma model defined as maps from the extended complex plane to rank-1 Hermitian projectors, with finite action, one can create a finite set of mutually orthogonal rank-1 Hermetian projectors  which are themselves solutions to the Euler-Lagrange equations \eqref{EL} with finite action \eqref{action}. We have also proven that the set of projectors contains an element that maps onto an equivalence class in $\CP$ which has a holomorphic representative and also one which has a anti-holomorphic representative.  Again, while there are analogous results for the vector form of the sigma model, the proofs and results of this paper have the benefit of fully utilizing the projective structures native to $\CP$ to both simplify the proofs and obtain new results as outlined in the Introduction.  

Next, we applied the proposed techniques and results to analyze surfaces created using the generalized Weierstrauss formula for immersion and were able to write several geometric characteristics of the surfaces, such as the metric tensor, area, Gaussian curvature and mean curvature in terms of physical quantities of the sigma model. Additionally, we gave a necessary and sufficient condition for surfaces to be associated with a $\CP$ model. Even from the two simple examples given in this paper, we can see that the types of surfaces obtained by this method are varied, having both positive and negative curvature, but still always conformally parameterized and with a finite area. The rich yet restrictive character makes such surfaces a rather special and interesting object of study. 

Furthermore, the results from this paper allow one to determine whether an arbitrary surface immersed in the $su(N)$ algebra is associated with a $\CP$ model and to derive physical quanties of the model from the geometry of the surface.  We anticipate that these results may be usefully applied to models of surfaces appearing in physics, chemistry and biology (see e.g. \cite{PolStrom1991,KonLand1999, CarKon1996, GPWbook, NPWbook, Saframbook, GinZinnbook, Raj2002, Davbook, Ou-YLuiXie1999,  Landolfi2003}), particularly in the cases where the associated analytical models are not fully developed. We hope to use these methods to shed some light on the matter. 

\end{section}
\begin{section}{Acknowledgements} This work was supported in part by research grants from NSERC of Canada. S.P. acknowledges a postdoctoral ISM fellowship awarded by the Mathematical Physics Laboratory of the Centre de Recherches Math\'ematiques, Universit\'e de Montr\'eal. This project was completed during A.M.G.'s visit to the \'Ecole Normale Sup\'erieure de Cachan (Centre de Math\'ematiques et leur Applications) and he would like to thank the CMLA for their kind invitation. 
\end{section}

\appendix
\section{Isometry between $\CP$ and the set of rank-1 Hermitian projectors}\label{AA}
\begin{property}
The set of vectors $f\in \mathbb{C}\Pb ^{N-1}$ are in bijection with the set of non-zero rank-1 projectors $\Pb .$ 
\end{property}
Proof. We give explicitly the map between the spaces. Given a non-zero rank-1 projector, there exists an arbitrary vector $\e$ such that $\Pb  \e\ne 0.$ The map and its inverse is given by 
$$\Pb =\Phi(f)=\frac{f\otimes f^\dagger}{f^\dagger f}, \qquad f=\Phi^{-1}(\Pb )=\Pb  e$$
We can see that this is well-defined for any vector in the equivalence class of $f$ since two equivalent vectors will have the same image, 
$$\vec{g}\sim \f \Rightarrow  \Phi(\f)=\frac{\f\otimes \f^\dagger}{\f^\dagger \f}=\frac{\vec{g}\otimes \vec{g}^\dagger}{\vec{g}^\dagger \vec{g}}=\Phi(\vec{g}).$$
Also, $\Phi^{-1}(\Pb )$ is independent of the choice of $\e:$ that is, for any $\v\in \C$ with $\Pb  \v\ne 0,$ $\Pb  \v $ is proportional to  $\Pb  \e$ and hence
$$\f\sim \Phi^{-1}(\Pb )\sim \Pb  \e \sim \Pb  \v \sim \f.$$
The two maps are mutual inverses 
$$\f\sim\Phi^{-1}(\Phi(\f))\sim \Phi^{-1}(\frac{\f\otimes \f^\dagger}{\f^\dagger \f})\sim \frac{\f^\dagger \e}{\f^\dagger \f}\f\sim \f.$$ 
$$\Pb =\Phi(\Phi^{-1}(\Pb ))=\Phi(\Pb  \e)=\frac{\Pb  \e\otimes \e^\dagger \Pb }{\e^\dagger \Pb  \e}=\Pb .$$
Thus the map between the two spaces is a bijection. \QED

\section{Proof of properties of rank-1 Hermitian projectors (1.1-1.3) }\label{AB}
 From any rank-1 Hermitian projector $\Pb $ we have the following identities
\begin{itemize}
 \item[1.1]  From any vector $\e,$ we can decompose $\Pb $ as proportional to the tensor product of $\Pb  \e$
\be \label{P11A} \Pb  \e\otimes \e^\dagger \Pb =\e^\dagger \Pb  \e\Pb \ee
If $\Pb  \e$ is also assumed to be nonzero, we can solve for $\Pb $ and obtain,
\be \label{PdecompA} \Pb =\frac{\Pb  \e\otimes \e^\dagger \Pb }{\e^\dagger \Pb  \e}.\ee 
\item[1.2] From any matrix $A,$ we have the following
\be \label{PtraceA}  \Pb  A\Pb =tr(\Pb  A)\Pb \ee
\item[1.3] Given, a mapping $\Pb \equiv \Pb (\xi_+, \xi_-): \mathbb{R}^2\rightarrow Gl(N, \mathbb{C})$ with $\Pb ^2=\Pb $ we have 
\be \label{dPA} \partial_\mu \Pb =\partial_{\mu}\Pb  \Pb +\Pb \partial_{\mu}\Pb , \qquad \partial_\mu=\partial_{\xi_\mu} \ee
\end{itemize}

Proof of 1.1  Assume $\Pb ^2=\Pb , \ \Pb ^\dagger=\Pb , $ and $tr(\Pb )=1.$ From any $\v\in \C, \ \Pb  \v=\Pb (\Pb  \v)$ and so $\Pb  \v$ is an eigenvector with eigenvalue 1. On the other hand, for any eigenvector $\w,$ we can use $\Pb ^2=\Pb $ to show
 \[  \Pb ^2\w=\Pb (\lambda \w)=\lambda^2 \w, \qquad \Pb ^2\w=\Pb  \w=\lambda \w.\]
But also, since $\Pb $ is Hermitian, its eigenvalues must be real and positive so $\lambda=1.$ Finally, since $tr(\Pb )=1,$ there has to be a unique normalized eigenvector. Thus, $\Pb  \e$ is either $0$ or in fact the unique eigenvector. In either case, we see that the following identity holds 
\[ \Pb  \e\otimes (\Pb  \e)^\dagger =(\e\Pb )^\dagger \Pb  \e\Pb .\]
\QED

Proof of 1.2 We use the decomposition above \eqref{Pdecomp} to expand
\[ \Pb  {\bf A} \Pb =\frac{\Pb  \e\otimes \e^\dagger \Pb  {\bf A}\Pb  \e\otimes \e^\dagger \Pb }{(\e^\dagger \Pb  \e)^2}=\frac{\e^\dagger \Pb  {\bf A}\Pb  \e}{\e^\dagger \Pb  \e}\Pb \]
Next we use the fact that 
\[\e^\dagger \Pb  {\bf A}\Pb  \e=tr(\Pb  {\bf A}\Pb  \e\otimes \e^\dagger)=tr({\bf A}\Pb )\e^\dagger \Pb  \e\]
where the final equality comes from the cyclic property of trace and \eqref{P11A}.
Thus, we see that $ \Pb  {\bf A}\Pb =tr(\Pb  {\bf A})\Pb .$ \QED
Proof of 1.3. This follows directly from differentiating $\Pb ^2=\Pb .$\QED
\section{Proof of properties of the sesquilinear product and the semi-norm }\label{AC}
Here we prove that if $\X$ and $\Y$ are arbitrary rank-1 Hermitian projectors and $\X \e, \Y \e \ne 0$, then 
\be\label{XYA}  \langle \X,\Y\rangle _{(\Pb )}=0 \iff \X \Y=\0.\ee

Proof. $\X^\dagger \Y=\0 \Rightarrow \langle \X,\Y\rangle =0$ follows from definition. On the other hand, assume $\langle \X,\Y\rangle =0.$
Since, $\X$ and $\Y$ are rank one projectors, $\X \v, \Y \v$ are proportional to $X \e, \Y \e$ respectively, for any vectors $\v,\w \in \C.$ But also
$$0=\langle \X,\Y\rangle =\frac{1}{\e^\dagger \Pb  \e}(\X \e)^\dagger(\Y \e) \propto (\X \v)^\dagger (\Y \w).$$ 
Thus $\v^\dagger \X\Y \w=0$ for any $\v,\w$ and so $\X\Y=\0.$ \QED  
\begin{property}\label{cauchyschwarz} The Cauchy-Schwarz inequality for the semi-norm $||\ ||_{(\Pb )}$ is given by 
$$|\langle \X,\Y\rangle _{(\Pb )}|\leq ||\X||_{(\Pb )}||\Y||_{(\Pb )}.$$
\end{property}
Proof,
we use the Cauchy-Schwarz inequality for the Euclidean norm 
$$|(\X \e)^\dagger(\Y \e)|\leq \sqrt{(\X \e)^\dagger \X \e}\sqrt{(\Y \e)^\dagger \Y \e}$$
to show
$$|\langle \X,\Y\rangle _{(\Pb )}|=|\frac{(\X \e)^\dagger(\Y \e)}{\e^\dagger \Pb  \e}|\leq\sqrt{ \frac{(\X \e)^\dagger \X \e}{\e^\dagger \Pb  \e}}\sqrt{\frac{(\Y \e)^\dagger \Y \e}{\e^\dagger \Pb  \e}}=||\X||_{(\Pb )}||\Y||_{(\Pb )}.$$ \QED

\begin{property}\label{triangle2} The triangle inequality for the semi-norm $||\ ||_{(\Pb )}$ is given by 
 $$||\X+\Y||_{(\Pb )}\leq ||\X||_{(\Pb )}+||\Y||_{(\Pb )}.$$
\end{property}
Proof, 
To prove this, we use the definition of the semi-norm, its sesquilinearity and the Cauchy-Schwarz inequality above to show, 
\bea ||\X+\Y||_{(\Pb )}^2&=&\langle \X,\X\rangle _{(\Pb )}+\langle \X,\Y\rangle _{(\Pb )}+\langle \Y,\X\rangle _{(\Pb )}+\langle \Y,\Y\rangle _{(\Pb )}\nn
&=&||\X||_{(\Pb )}^2+||\Y||_{(\Pb )}^2+2Re(\langle \X,\Y\rangle _{(\Pb )})\nn
&\leq &||\X||_{(\Pb )}^2+||\Y||_{(\Pb )}^2+2|\langle \X,\Y\rangle_{(\Pb )}|\nn
&\leq & ||\X||_{(\Pb )}^2+||\Y||_{(\Pb )}^2+2||\X||_{(\Pb )}||\Y||_{(\Pb )}=(||\X||_{(\Pb )}+||\Y||_{(\Pb )})^2.\nonumber \eea \QED

\section{Proof of properties of the covariant derivative (4.1-4.5)}\label{AD}
In this section, in order to lighten the notation we have fixed $P$ a rank-1 Hermitian projector which is used to define the sesquilinear product, semi-norm and covariant derivative. We drop the subscripts and superscripts denoting this dependence as in the body of the paper. 
For the given rank-1 Hermitian projector $\Pb $ and arbitrary matrix functions $\X,\Y$ the following hold:
\begin{itemize}
\item[4.1]  $\Pb $ and  $D_\pm \Pb $ are orthogonal\\
\be \label{PDPA} \langle \Pb ,D_{\pm}\Pb \rangle=0.\ee
\item [4.2] The covariant derivative is compatible with the sesquilinear product in the following sense
 \be \label{d<>A}\frac{\partial}{\partial \mu} \langle \X,\Y\rangle =\langle D_{\overline{\vec{\mu}}}\X,\Y\rangle+\langle \X,D_{\vec{\mu}}\Y\rangle\ee
 \item[4.3] If $\X$ is a rank-1 Hermitian projector,  there exists a scalar function $\phi(\X)$ so that 
\be \label{DXA} \partial_{\pm} \X \X \e=(D_{\pm}\X-\phi(\X) \X)\e\ee
 and in particular, if $\X=\Pb $ we have 
 \be \label{DPA} \partial_{\pm} \Pb  \Pb  \e=D_{\pm}\Pb  \e.\ee
 
 \item[4.4]The commutator of covariant derivatives is given by
$$[D_{+},D_{-}]\X=(||D_+\Pb ||^2-||D_-\Pb ||^2) \X.$$
\item[4.5] The norm of $D_\pm \Pb $ is given by
 \be \label{normDPA}|| D_+\Pb ||^2=tr(\p \Pb  \Pb  \bp \Pb ) \qquad ||D_-\Pb ||^2=tr(\bp \Pb  \Pb  \p)\ee
 \end{itemize}

Proof of 4.1. This follows from direct computation. 
$$\langle \Pb ,D_{\pm }\Pb \rangle =\frac{1}{\e^\dagger \Pb  \e}\left(\e^\dagger \Pb  \partial_\pm \Pb  \e-\frac{\e^\dagger \Pb \partial_\pm \Pb  \e}{\e^\dagger \Pb \e}\e^\dagger \Pb  \e\right)=0$$
\QED
Proof of 4.2. We compute
\bea\partial_\mu\langle \X,\Y\rangle&=&\partial_\mu\left(\frac{\e^\dagger \X^\dagger \Y \e}{\e^\dagger \Pb  \e}\right)\nn
&=&\frac{\e^\dagger \left(\partial_{\overline{\mu}}\X\right)^\dagger \Y \e}{\e^\dagger \Pb  \e}+\frac{ \e^\dagger \X^\dagger \partial_\mu \Y  \e}{\e^\dagger \Pb  \e}-\frac{\e^\dagger \X^\dagger \Y \e}{(\e^\dagger \Pb  \e)^2}\e^\dagger\partial_\mu \Pb  \e \nonumber .\eea
Next, we use the fact that $\partial_\mu \Pb =\Pb \partial_\mu \Pb +(\partial_\mu \Pb )\Pb $
to show that 
\bea\frac{\e^\dagger \X^\dagger \Y \e}{(\e^\dagger \Pb  \e)^2}\e^\dagger\frac{\partial}{\partial \mu}\Pb  \e&=&
\frac{\e^\dagger \left((\e^\dagger \Pb \partial_{\overline{\mu}}\Pb  \e)\X\right)^\dagger \Y \e}{(\e^\dagger \Pb  \e)^2}+\frac{\e^\dagger \X^\dagger (\e^\dagger \Pb \partial_\mu \Pb  \e) \Y \e }{(\e^\dagger \Pb  \e)^2}\nonumber 
\eea
and so we obtain
\bea \partial_\mu \langle \X,\Y\rangle&=& \frac{\e^\dagger \left(\partial_{\overline{\mu}}\X-\langle \Pb ,\partial_{\overline{\mu}}\Pb  \rangle \X\right)^\dagger \Y \e}{\e^\dagger \Pb  \e}+\frac{ \e^\dagger \X^\dagger \left(\partial_\mu \Y-\langle \Pb ,\partial_\mu \Pb  \rangle \Y \right) \e}{\e^\dagger \Pb  \e} \nn
&=&\langle D_{\overline{\vec{\mu}}}\X,\Y\rangle+\langle \X,D_{\vec{\mu}}\Y\rangle. 
\nonumber\eea \QED

Proof of 4.3 . We use  \eqref{DPA} to commute $\X$ with $\partial_\pm \X$ and the definition of the covariant derivative to obtain
\be\label{1A} \partial_{\pm} \X \X \e=(\partial_{\pm} \X-\X\partial_{\pm} \X)\e=(D_{\pm}\X+\langle \Pb , \partial_\pm \Pb \rangle \X-\X\partial_{\pm} \X)\e.\ee
Next, we use  \eqref{P11A} to expand $\X$ as 
\[\X=\X \e\otimes \e^\dagger \X /(\e^\dagger \X \e).\]
 Substituting this into the right hand side of \eqref{1A}, we obtain
\bea \label{2A} \partial_{\pm} \X \X \e&=&(D_{\pm}\X+\langle \Pb , \partial_\pm \Pb \rangle \X)\e-\frac{\X \e\otimes \e^\dagger \X}{\e^\dagger \X \e}\partial_{\pm} \X \e\nn
\label{3A}&=&(D_{\pm}\X+\langle \Pb , \partial_\pm \Pb \rangle \X+\frac{\e^\dagger \X\partial_{\pm} \X \e}{\e^\dagger \X \e}\X)\e.\eea
As in the equation above, we define the scalar function 
$$\phi(\X)=\langle \Pb , \partial_\pm \Pb \rangle+\frac{\e^\dagger \X\partial_{\pm} \X \e}{\e^\dagger \X \e},$$
and then \eqref{3A} becomes 
\be \partial_{\pm} \X \X \e=(D_{\pm}\X+\phi(\X) \X)\e.\ee 
If we consider the case that $\X=\Pb $ we compute
$$\phi(\Pb )=\langle \Pb , \partial_\pm \Pb \rangle+\frac{\e^\dagger \Pb \partial_{\pm} \Pb  \e}{\e^\dagger \Pb  \e}=0.$$
 \QED 
 Proof of 4.4. From the definition of covariant derivative, we obtain
\[ [D_{+},D_{-}]\X=(\bp \langle \Pb , \p \Pb  \rangle -\p \langle \Pb , \bp \Pb  \rangle) = ( \langle D_+ \Pb , \p \Pb  \rangle - \langle D_- \Pb , \bp \Pb  \rangle) \]
but also, $\partial _\pm \Pb  =D_\pm \Pb  -\langle \Pb  \partial_\pm \Pb \rangle \Pb $ and then we use the fact that $\Pb $ is orthogonal to $D_\pm \Pb $ to obtain $\langle D_\pm \Pb , \partial_\pm \Pb \rangle =\langle D_\pm \Pb , D_\pm \Pb \rangle$ and so 
$$[D_{+},D_{-}]\X=( \langle D_+\Pb  , D_+\Pb  \rangle +\langle D_- \Pb , D_- \Pb \rangle )\X= (||D_+\Pb ||^2-||D_-\Pb ||^2) \X.$$
 \QED
Proof of 4.5 . We use  \eqref{DPA} and   \eqref{PtraceA}  which gives $\Pb  A\Pb =tr(\Pb  A)\Pb ,$ to obtain
\bea ||D_+\Pb ||^2&=&\frac{(D_+\Pb  \e)^\dagger D_+\Pb  \e}{\e^\dagger \Pb  \e}\nn
&=& \frac{\e^\dagger \Pb  \bp \Pb  \p \Pb  \e}{\e^\dagger \Pb  \e} = tr(\Pb  \bp \Pb  \p \Pb ).\eea
The computation for 
\be ||D_-\Pb ||^2=tr(\Pb  \p \Pb  \bp \Pb )\ee
 is identical. We use the cyclic property of trace to rewrite the quantities as 
\be || D_+\Pb ||^2=tr(\p \Pb  \Pb  \bp \Pb ), \qquad ||D_-\Pb ||^2=tr(\bp \Pb  \Pb  \p \Pb).\ee \QED
\section{Proof of holomorphicity properties (5.1-5.3)}\label{AE}

The following properties of holomorphic vectors and projectors hold.
\begin{itemize}
\item[5.1] A vector $\x$ in  $\CP$ is holomorphic if and only if $\bp \x$ is proportional to $ \x.$
\item[5.2] A Hermitian, rank one projector $\X$ is holomorphic if and only if 
\be \label{projholoA} \X\partial_+ \X=\partial_-\X\X=\0. \ee 
\item[5.3] If $\X=\Pb $ the fixed projector used in the definition of the covariant derivative \eqref{D}, then $\Pb $ is holomorphic if and only if 
\be D_-\Pb \e =\0.\ee
\end{itemize}

Proof of 5.1
$\x$ is holomorphic if and only if there is some non-zero scalar function $c(\xi_{+},\xi_{-}): \mathbb{R}^2\rightarrow \mathbb{C}$ for which 
$$\bp (c\x)=\bp(c) \x+c\bp \x=\0 \iff  \bp \ x=\frac{\bp c}{c} \x.$$
\QED
Proof of 5.2. Given an arbitrary rank-1 Hermitian projector $\X$, let $\vec{x} \in \CP(\xi_+, \xi_-)$ be the non-zero vector on to which $\X$ projects. That is, $\X \vec{v}$ is proportional to $ \x$ for any vector $\vec{v}.$ If $\vec{x}$ is holomorphic, then $\bp(c\vec{x})=0.$ Define $\hat{\vec{x}}=c\vec{x}, $ then 
$$\X=\frac{\hat{\vec{x}}\otimes \hat{\vec{x}}^\dagger}{\hat{\vec{x}}^\dagger \hat{\vec{x}}},\qquad \partial_-\hat{\vec{x}}=\0.$$ By direct computation we see that $\X\partial_+ \X=\0.$
On the other hand, if we assume that $\X\partial_+ \X=\0$ we obtain
$$\0=\X\partial_+ \X=\frac{\vec{x}\otimes(\bp \vec{x}-\frac{(\bp \vec{x})^\dagger \vec{x}}{\vec{x}^\dagger \vec{x}})\vec{x}}{\vec{x}^\dagger \vec{x}}$$
which implies
$$ \bp \vec{x}-\frac{(\bp \vec{x})^\dagger \vec{x}}{\vec{x}^\dagger \vec{x}}\vec{x}=\0.$$
\QED
Proof of 5.3. If $\X=\Pb $ then $D_-\Pb  \e=\bp \Pb  \Pb  \e$ which equals $\0$ by 4.2. \QED
The analogous statements for anti-holomorphicity follow directly. That is, $\X$ is antiholomorphic if the projector which it projects onto is antiholomorphic and this is equivalent to 
\be \X\partial_- \X=\partial_+\X\X=\0.\ee

\section{Proof of properties of the raising and lowering operators (6.1-6.4)} 
\label{AF}

From $\X$ an Hermitian rank-1 projector function the following hold. 
\begin{itemize}
 \item[6.1] $\Pi_{\pm}\X$ is a well defined Hermitian projector which projects onto the vector $ \partial_\pm \X \X e.$
 Thus, either $\Pi_{\pm}\X$ is $\0$ or it is of rank-1 and can be written as 
\be \label{PXA} \Pi_{\pm}\X=\frac{ \partial_\pm \X \X \e\otimes \e^\dagger \X \partial_\mp \X}{\e^\dagger \X\partial_\mp \X\partial_\pm \X \X \e}.\ee
\item[6.2]  $\X$ is orthogonal to $\Pi_{\pm}\X.$
\item[6.3] The raising and lowering operators acting on $\Pb $ can be written in terms of the covariant derivative 
\be\label{rlcovA} \Pi_{\pm}\Pb =\frac{D_{\pm}\Pb  \e\otimes \e^\dagger (D_{\pm}\Pb )^\dagger}{(D_{\pm}\Pb  \e)^\dagger D_{\pm}\Pb  \e}\ee
\item[6.4]  $\X$ is holomorphic if and only if 
\be \Pi_{-}\X=\0.\ee
\end{itemize}
Proof of 6.1.
First, we note that since $X$ is a Hermitian projector $tr(\partial_\pm \X \X \partial_\mp \X)$ can be written as $tr({\bf A}^\dagger {\bf A})$ with ${\bf A}= \X \partial_\mp \X$ and so will be identically $0$ only when $\X \partial_\mp \X=\0, $ thus the raising and lowering operators are well defined when they act Hermitian projectors. \\
Next, assume $\Pi_\pm \X =\0.$ We see immediately that it is Hermitian and 
\[ \partial_\pm \X \X \e \otimes \e^\dagger \X \partial_\mp \X =\0 \]
which implies that $\partial_\pm \X \X e=\0$ and so $\Pi_\pm \X =\0$ maps onto $\partial_\pm \X \X e=\0.$
Assume now that $\Pi_\pm \X \ne \0.$ Using \eqref{P11A} and \eqref{PtraceA} one finds $(\Pi_{\pm}\X)^2=\Pi_{\pm}\X, \ (\Pi_{\pm}\X)^\dagger=\Pi_{\pm}\X$ and for any vector $\v$
$$\Pi_{\pm}\X \v =\frac{\e^\dagger \X \partial_{\mp} \X \v}{\e^\dagger \X\bp \X\partial \X \X \e} \partial_{\pm} \X \X \e.$$
Thus, $\Pi_\pm \X$ projects onto $\partial_\pm  \X \X \e$ can so can be written as \eqref{PXA}. \QED
Proof of 6.2. This follows directly from differentiating the projective property of $\X.$
\[ \partial_\pm \X=\X \partial_\pm \X+\partial_\pm \X \X\Rightarrow \X\partial_\pm \X \X=\0\]
and so 
\[ \X\Pi_{\pm}\X =\frac{\X\partial_\pm \X \X\partial_\mp \X}{tr(\partial_\pm \X \X\partial_\mp \X)}=\0.\] \QED
Proof of 6.3. This follows immediately from  \eqref{PXA} and  \eqref{DPA}  which relate the covariant derivative of $\Pb $ with $\Pb \partial _\pm \Pb $. \QED

Proof of 6.4. It has been proven in \cite{GoldGrund2010} that the lowering operator annihilates holomorphic projectors and here we also prove the contrapositive.  Given a Hermitian, rank-1 projector $\X. $ If it is holomorphic then by  \eqref{projholoA} 
$$\Pi_{-}\X=\frac{\bp \X \X\partial_+ \X}{tr(\bp \X \X\partial_+ \X)}=\0.$$
From the other direction, suppose that $\Pi_{-}\X=\0$ then
$$(\X\partial \X)^\dagger \X\partial \X=\0 \Rightarrow \X\partial \X=\0 $$
and so $\X$ is holomorphic. 
\QED

\section{Proof that raising and lowering operators are mutual inverses on solutions of the Euler-Lagrange equations}\label{AG}
We prove that if $\Pb $ satisfies Euler-Lagrange equation \eqref{EL}, then 
\[\Pi_\pm(\Pi_\mp \Pb )=\Pb ,\qquad \Pi_\mp(\Pi_\pm \Pb  )=\Pb, \] 
whenever the term inside the parentheses is not $\0.$
Proof. Assume $\Pi_+\Pb  \ne \0$, we will show $\Pi_-\Pi_+\Pb =\Pb .$ By definition of the raising operators
\[ \Pi_-\Pi_+\Pb =\frac{\bp \Pi_+\Pb  \Pi_+\Pb  \p \Pi_+\Pb }{tr(\bp \Pi_+\Pb  \Pi_+\Pb  \p \Pi_+\Pb )}\]
and so we rewrite  $\Pi_-\Pi_+\Pb =\X\X^\dagger$ with $\X=(tr(\bp \Pi_+\Pb  \Pi_+\Pb  \p \Pi_+\Pb ))^{-1/2}\bp \Pi_+\Pb  \Pi_+\Pb.$
Thus, if we  prove $\X=\Pb  A$ for some matrix ${\bf A},$ we will have 
$$ \Pi_-\Pi_+\Pb =\Pb {\bf A}{\bf A}^\dagger \Pb =tr(\Pb  {\bf A}{\bf A}^\dagger)\Pb $$
but $tr(\Pi_-\Pi_+\Pb )=1$ and so $tr(\Pb  {\bf A}{\bf A}^\dagger)=1$ and we will have proven the theorem. 
We compute 
\bea \bp \Pi_+\Pb \Pi_+\Pb  &=&\frac{\bp \p \Pb  \Pb  \bp \Pb +\p \Pb  \bp \Pb  \bp \Pb  +\p \Pb  \Pb  \bp^2\Pb }{tr(\p \Pb  \Pb \bp \Pb  )}\Pi_+\Pb \nn&-&\frac{tr(\bp \p \Pb  \Pb  \bp \Pb +\p \Pb  \bp \Pb  \bp \Pb  +\p \Pb  \Pb  \bp^2)\p \Pb  \Pb  \bp \Pb }{ tr(\p \Pb  \Pb \bp \Pb  )^2}\Pi_+\Pb \nonumber\eea
Since $[\bp \p \Pb  ,\Pb  ]=\0,$  $\bp \p \Pb  \Pb  \bp \Pb =\Pb  \bp \p \Pb  \bp \Pb $ and also $tr(\bp \p \Pb  \Pb  \bp \Pb )=0.$ Using  \eqref{dP}  we  compute
\[ tr(\p \Pb  \bp \Pb  \bp \Pb )=tr(\p \Pb  \Pb  \bp \Pb  \bp \Pb  \Pb ) +tr(\p \Pb  \bp \Pb  \Pb  \bp \Pb )=0.\]
We see that the folllowing terms cancel 
$$\left(\frac{\p \Pb  \Pb  \bp^2}{tr(\p \Pb  \Pb \bp \Pb  )}-\frac{tr(\p \Pb  \Pb  \bp^2\Pb )\bp \Pb  \Pb  \p \Pb }{ tr(\p \Pb  \Pb \bp \Pb  )^2}\right)\frac{\p \Pb  \Pb \bp \Pb }{tr(\p \Pb  \bp \Pb  )} =\0.$$
So that finally, we have 
\bea \bp \Pi_+\Pb \Pi_+\Pb &=&\left(\frac{\Pb \bp \p \Pb  \bp \Pb +\p \Pb  \bp \Pb  \bp \Pb }{tr(\p \Pb  \Pb \bp \Pb  )}\right)\Pi_+\Pb \nn
&=&\Pb \left(\frac{\bp \p \Pb  \bp \Pb \p \Pb  \bp \Pb  +tr(\Pb \p \Pb  \bp \Pb )tr(\Pb  \bp \Pb \p \Pb )\bp \Pb }{tr(\p \Pb  \Pb \bp \Pb  )^2}\right).\nonumber\eea 
Thus $\Pi_-\Pi_+\Pb =\Pb. $ By direct analogy we see that  $\Pi_+\Pi_-\Pb =\Pb $ as well. \QED
 \bibliography{sigmamodels}{}
\bibliographystyle{model1-num-names}

\end{document}